\documentclass{article}
\usepackage{amsthm,amsfonts, amsbsy, amssymb,amsmath,graphicx}
\usepackage{graphics}

\newtheorem{thm}{Theorem}
\newtheorem{lm}{Lemma}

\title{On two categorifications of the arrow polynomial for
virtual knots}

\author{Heather Ann Dye, Louis Hirsch Kauffman\footnote{Corresponding Author}, Vassily Olegovich Manturov
\footnote{Corresponding Author}}

\begin{document}

\maketitle

\newcommand{\skcrossr}{\raisebox{-0.25\height}{\includegraphics[width=0.5cm]{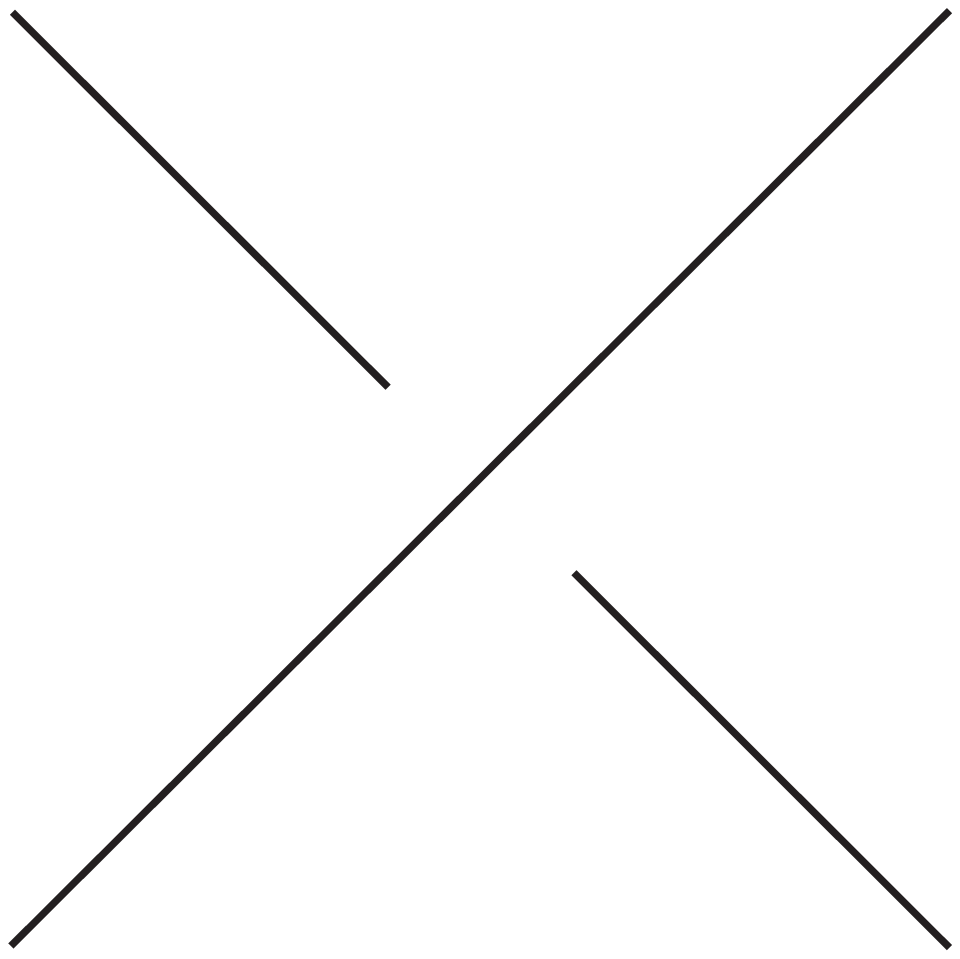}}}
\newcommand{\skkinkr}{\raisebox{-0.25\height}{\includegraphics[width=0.5cm]{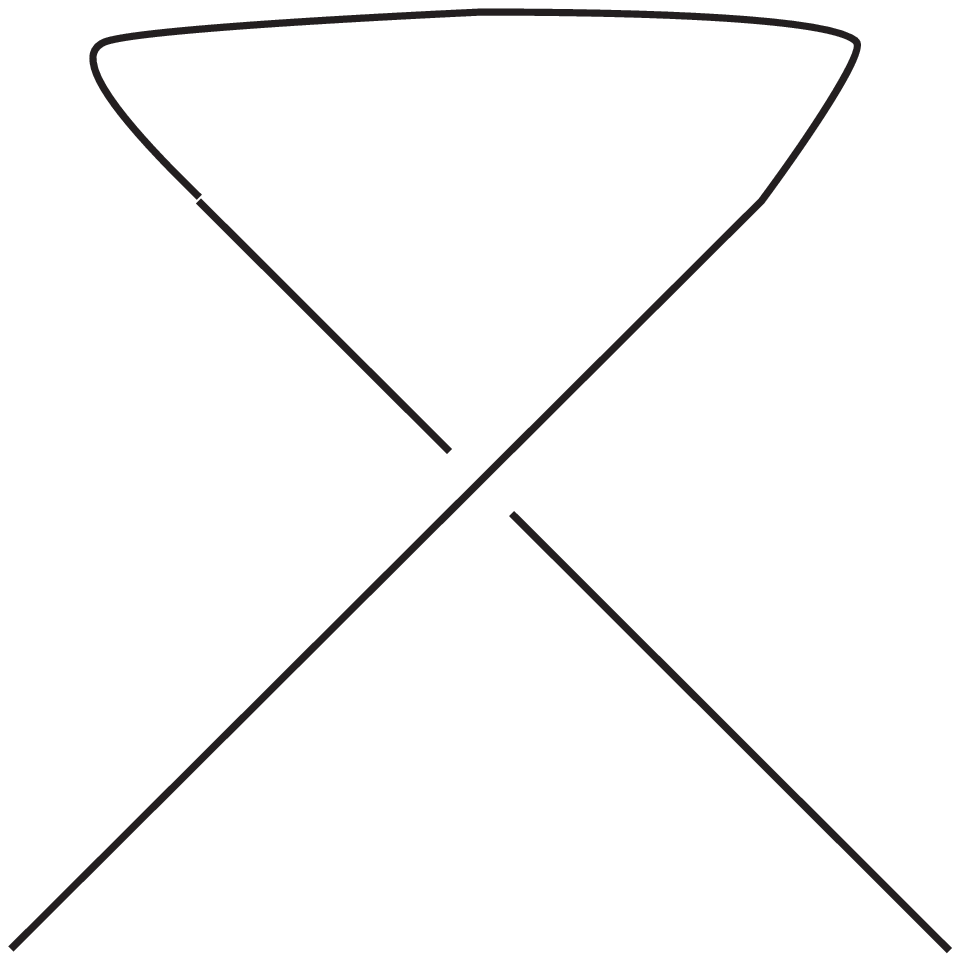}}}
\newcommand{\skkinkl}{\raisebox{-0.25\height}{\includegraphics[width=0.5cm]{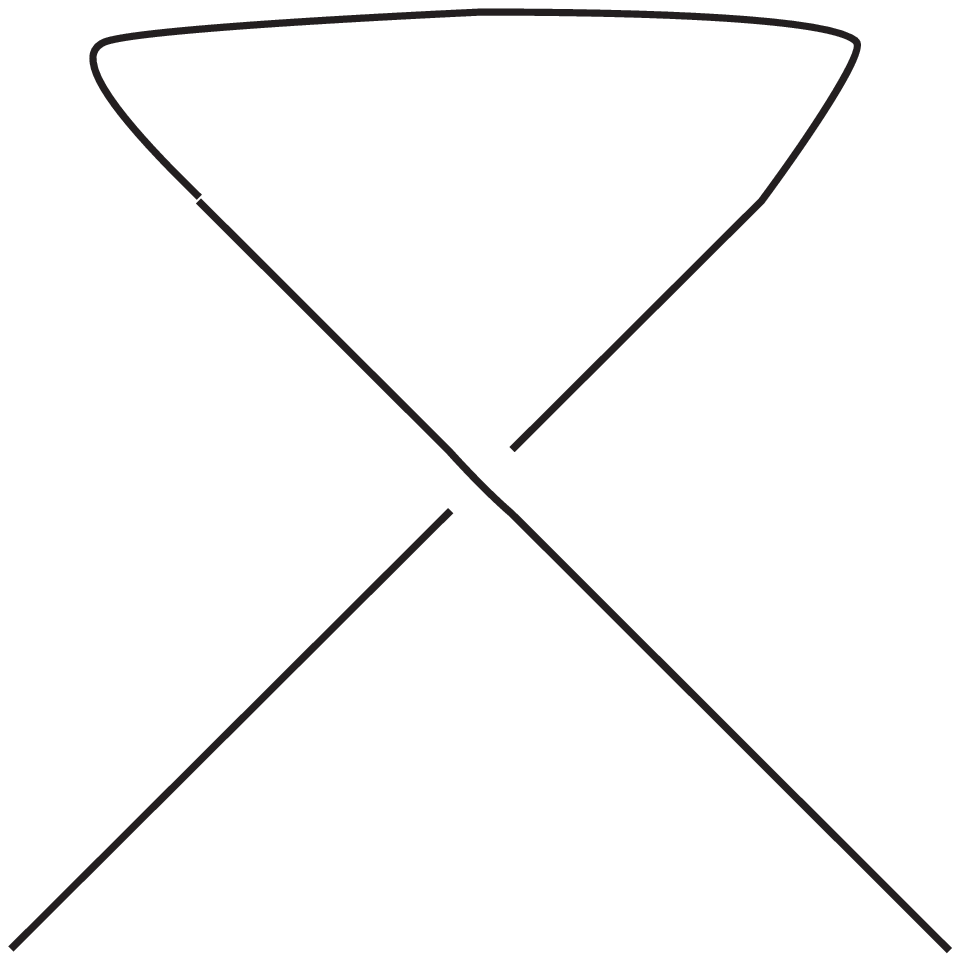}}}
\newcommand{\skroh}{\raisebox{-0.25\height}{\includegraphics[width=0.5cm]{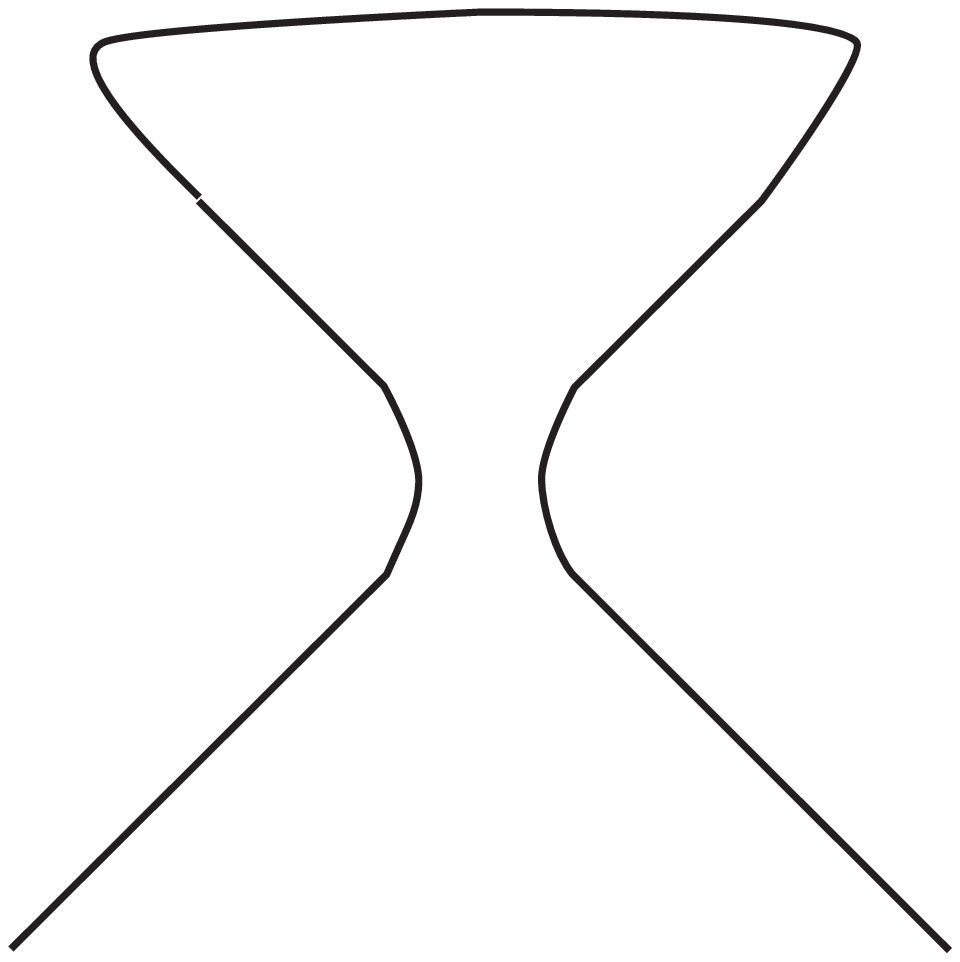}}}
\newcommand{\skrov}{\raisebox{-0.25\height}{\includegraphics[width=0.5cm]{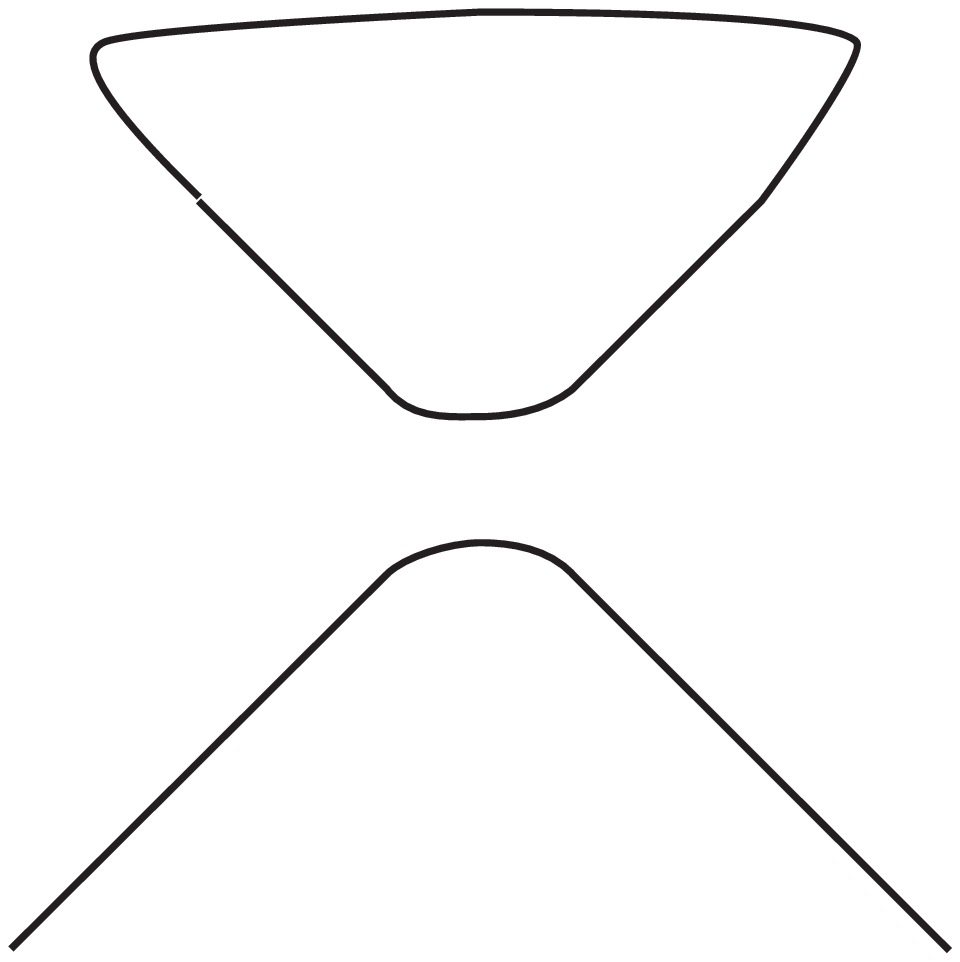}}}
\newcommand{\skrtwhh}{\raisebox{-0.25\height}{\includegraphics[width=0.5cm]{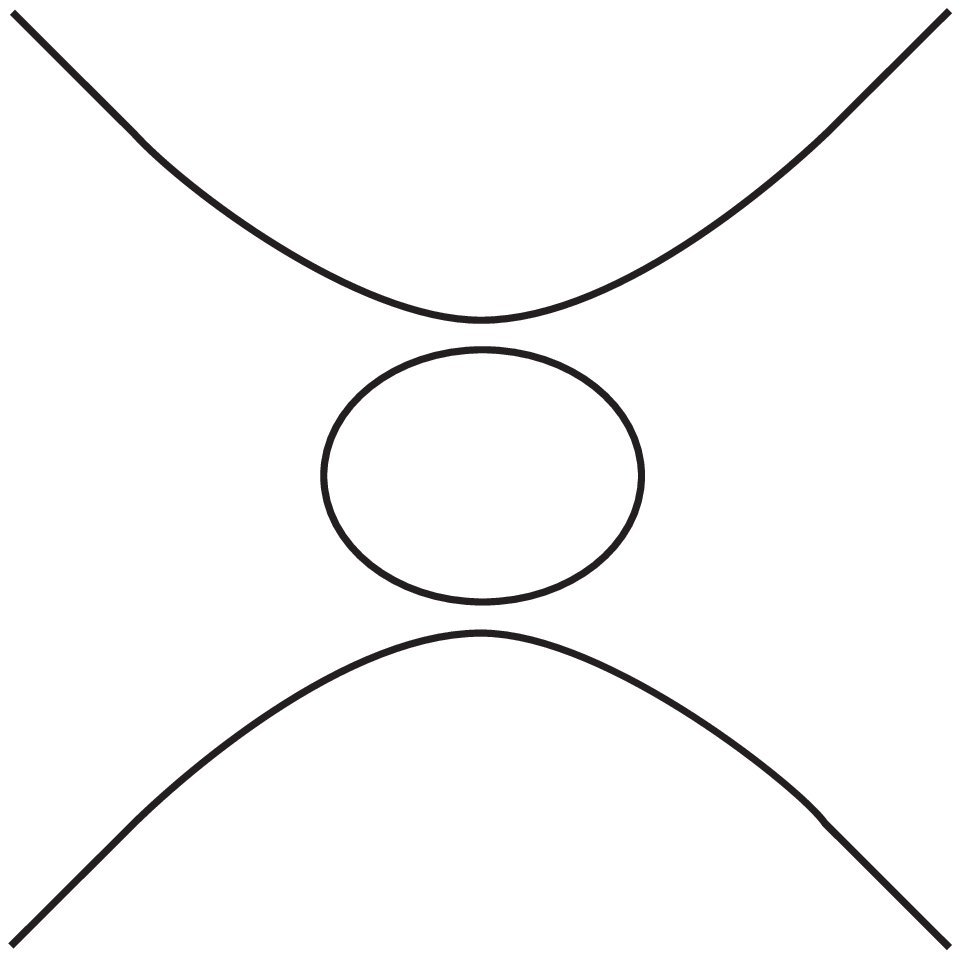}}}
\newcommand{\skrtwvh}{\raisebox{-0.25\height}{\includegraphics[width=0.5cm]{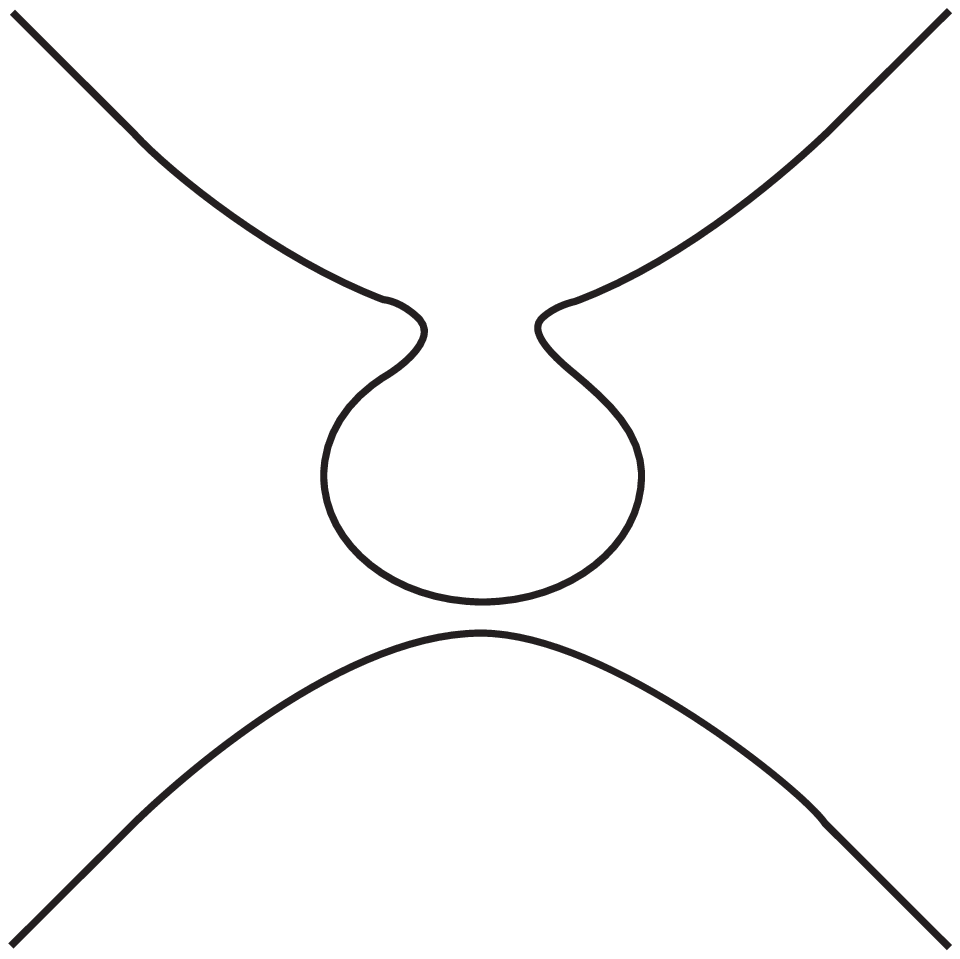}}}
\newcommand{\skrtwhv}{\raisebox{-0.25\height}{\includegraphics[width=0.5cm]{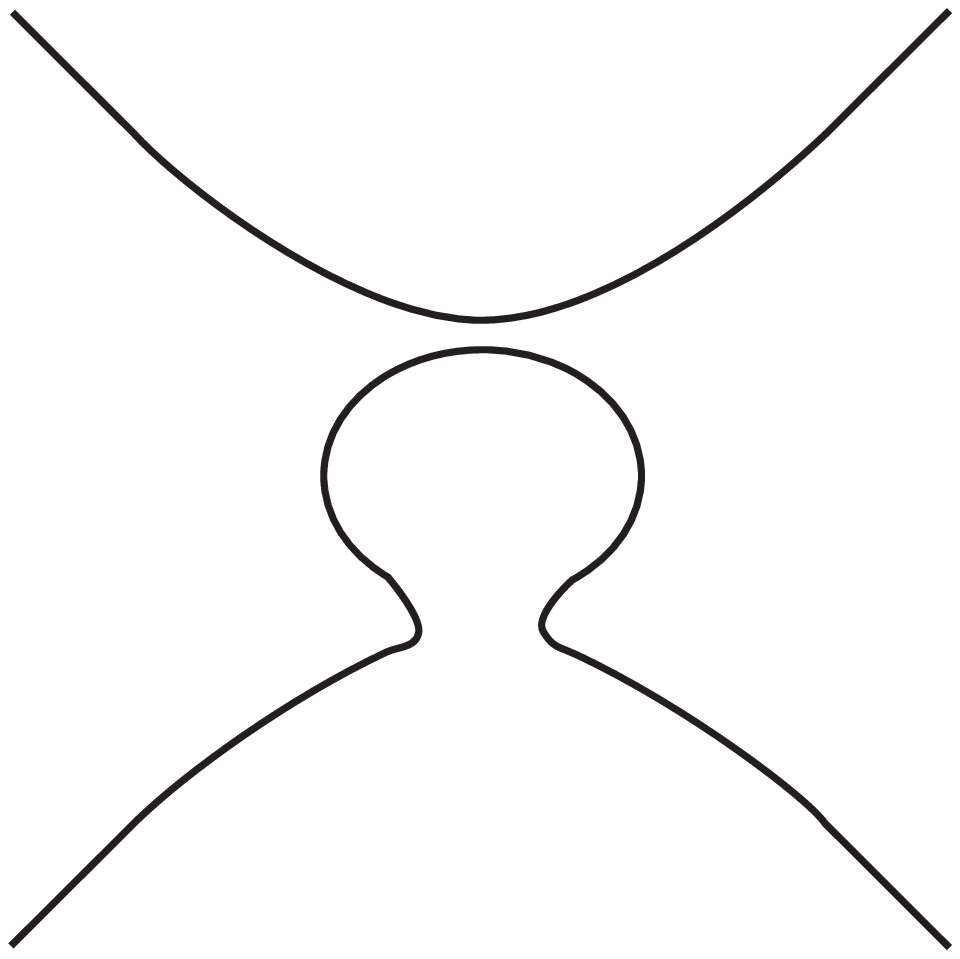}}}
\newcommand{\skrtwvv}{\raisebox{-0.25\height}{\includegraphics[width=0.5cm]{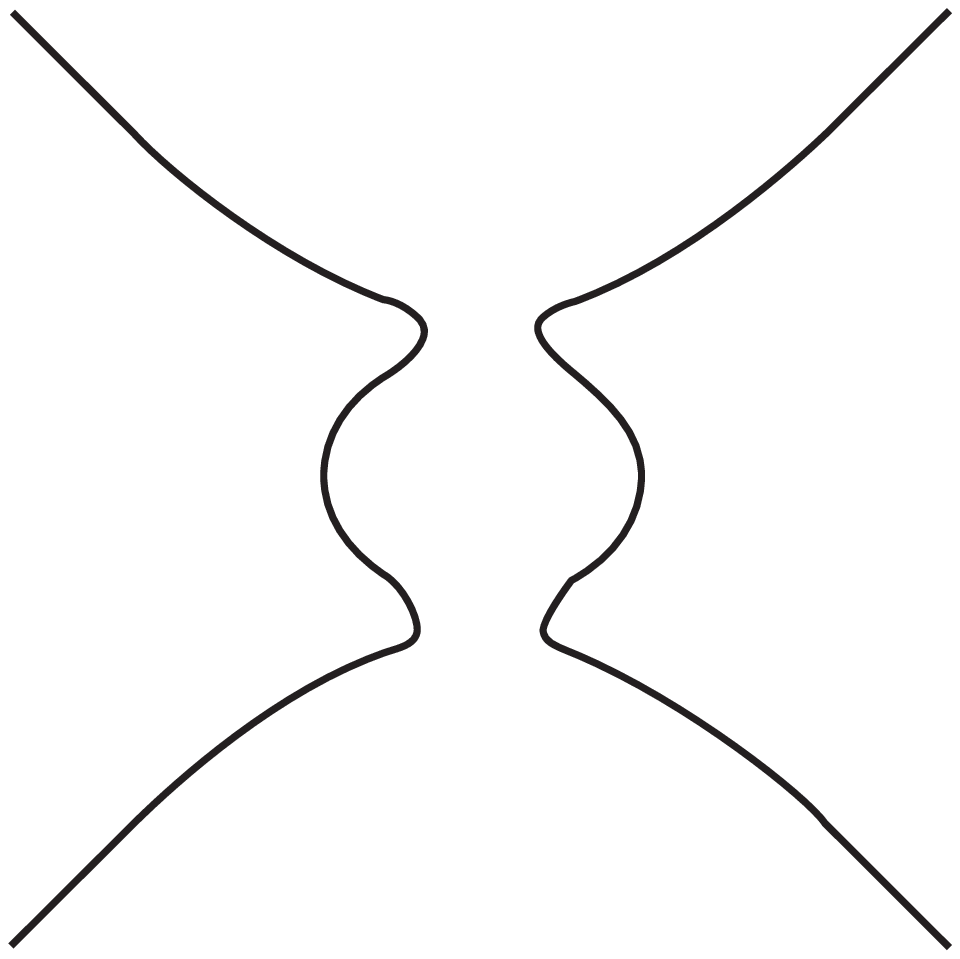}}}
\newcommand{\skcrv}{\raisebox{-0.25\height}{\includegraphics[width=0.5cm]{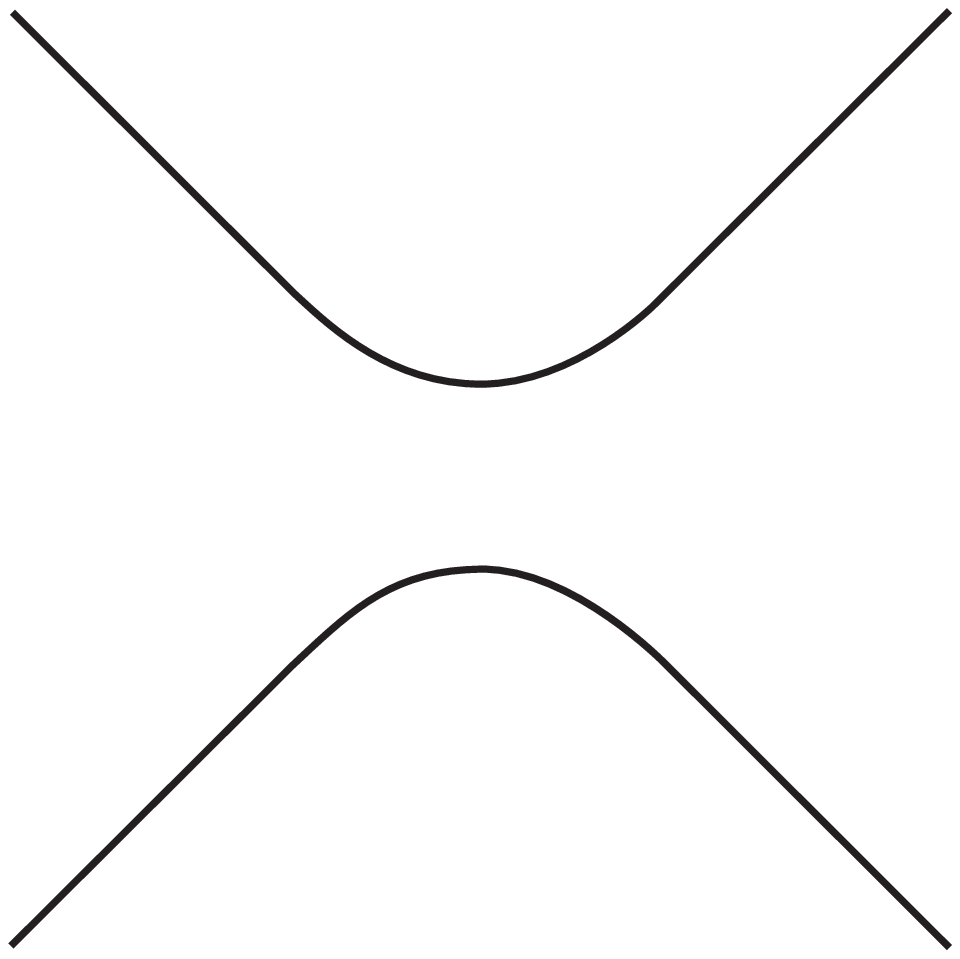}}}
\newcommand{\skcrh}{\raisebox{-0.25\height}{\includegraphics[width=0.5cm]{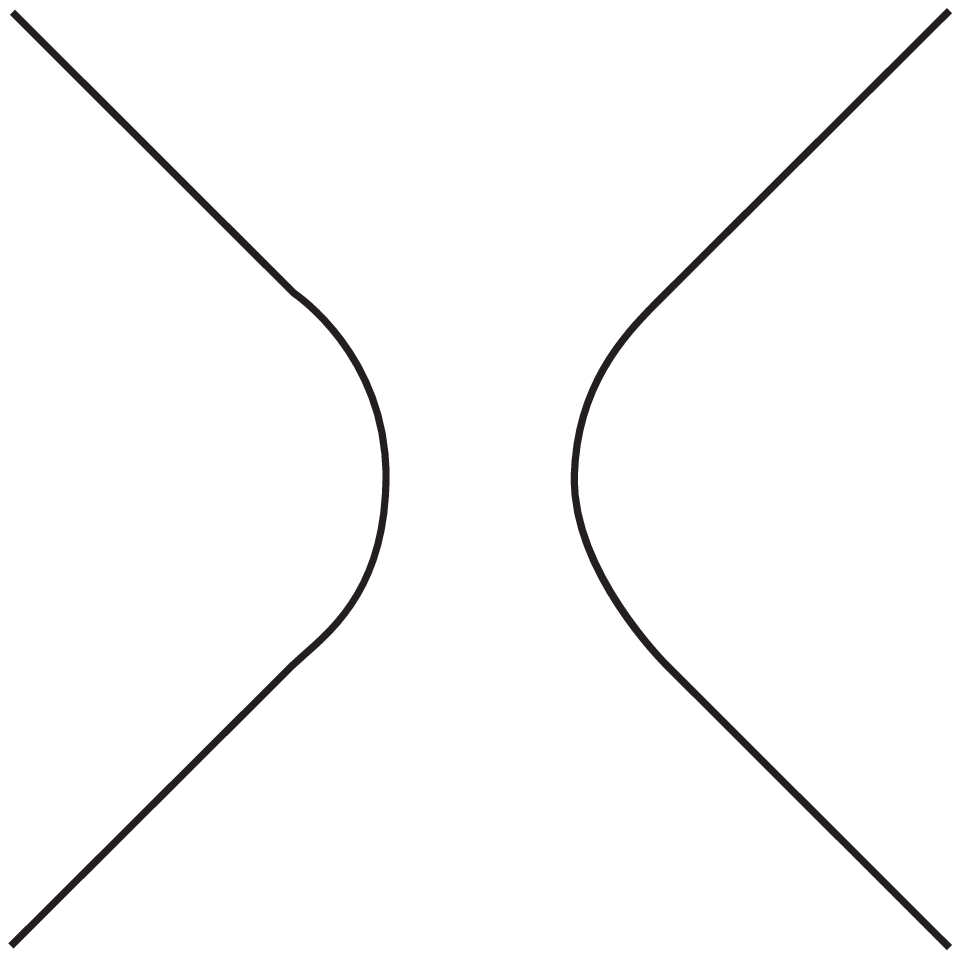}}}

\newcommand{\ddO}{\mbox{$\dot{1}$}}
\newcommand{\ddX}{\mbox{$\dot X$}}
\newcommand{\dg}{\mbox{deg}}
\newcommand{\lb}{\mbox{[[}}
\newcommand{\rb}{\mbox{]]}}

\abstract{
Two categorifications are given for the arrow polynomial, an
extension of the Kauffman bracket polynomial for virtual knots. The arrow
polynomial extends the bracket polynomial to infinitely many variables,
each variable corresponding to an integer {\it arrow number} calculated
from each loop in an oriented state summation for the bracket. The
categorifications are based on new gradings associated with these arrow
numbers, and give homology theories associated with oriented
virtual knots and links via extra structure on the Khovanov
chain complex. Applications are given to the estimation
of virtual crossing number and surface genus of virtual knots and links.

 }


\section{Introduction}
The purpose of this paper is to give a categorification for an
extension of the Kauffman bracket polynomial, giving a new
categorified homology for virtual knots and links. The extension of
the bracket that we work with is the {\it arrow polynomial} as
defined in \cite{ExtBr,DK}. This invariant was independently
constructed by Miyazawa in \cite{Miya,Miy} and so this work can also
be seen as a categorification of the Miyazawa polynomial. \bigbreak

In \cite{ExtBr}, Kauffman gives an extension of the bracket
polynomial for virtual knots that is obtained by using an oriented
state expansion, as indicated here in Figure \ref{Figure 1}. In such
an expansion there are two types of smoothing as shown in this
figure. The guiding principle for the extended bracket invariant is
to retain the pairing of the cusps at the reverse oriented
smoothings for as long as possible. The resulting state
configurations are then replaced by 4-regular virtual graphs, and
the invariant is a linear combination of these graphs with
polynomial coefficients. For a given state $S$, the corresponding
graph is denoted by $[S].$ In \cite{ExtBr} this invariant is then
simplified by retaining the cusps at the non-oriented smoothings but
not insisting upon pairing them. In this simplified version $[S]$ is
replaced by a diagram that is a union of circle graphs with
(reduced) cusps and virtual crossings, modulo virtual equivalence.
States are reduced via the rule that consecutive pairs of cusps on a
given state curve cancel if they point to the same local side of the
curve in the plane.  With this caveat, each state curve can be
regarded as an extra variable $K_{n}$ with an index $n$ denoting one
half of the reduced number of cusps. This simplified version of the
invariant is called the {\it arrow polynomial}. It takes the form

\begin{equation}
\langle K\rangle_{A} = \sum_{S}A^{\alpha(S)-\beta(S)}(-A^{2}-A^{-2})^{\gamma(S)-1} \prod_{c \in
c(S)}K_{n(c)},\label{secondbracket}
\end{equation}

\noindent where the product is taken over all single loop components $c$ in the state $S$, and
$n(c)$ counts one half the number of cusps in the reduced circle graph. Here $\alpha(S)$ and $\beta(S)$ are the numbers
of positively (resp., negatively) smoothed crossings, and
$\gamma(S)$ is the number of loops in the state $S.$
\bigbreak

H. A. Dye and L. H. Kauffman studied an equivalent version of
the arrow polynomial \cite{DK} and used it to obtain a
lower bound on the virtual crossing number for diagrams of a virtual link.
In the Dye-Kauffman version the cusps are replaced by an extra
orientation convention. See Figure \ref{Figure 3}. Here we shall refer to $ \langle K\rangle_{A}$ as the
arrow polynomial of $K.$  We call the reduced cusp count
$n(c)$ for a state loop the {\it arrow number} of this loop. Thus a state loop
with label $K_{n}$ has arrow number $n.$ \bigbreak

Both the extended bracket polynomial and the arrow polynomial $\langle K \rangle_{A}$ are invariant
with respect to the second and the third Reidemeister moves. They can be made invariant under
the first Reidemeister move by the usual normalisation by a power of
$(-A^{3}).$ In the rest of the paper, we shall omit this normalisation. Moreover,
while passing to the Khovanov homology, we shall omit the
corresponding renormalisation and refer the interested reader to
\cite{BN,MyBook} or \cite{BN,BN2}. \bigbreak

The aim of the present paper is to present two categorifications of
the arrow polynomial \cite{ExtBr,DK}. We split the chain spaces of
the Khovanov complex ${\cal C}(K)$ into subspaces ${\cal
C}_{gr=x}(K)$ with a fixed new grading $x$ and restrict our
differential $\partial$ to these subspaces. Now, set
$\partial=\partial'+\partial''$ where $\partial'$ is the part of
$\partial$ which preserves the new gradings for basic chains, and
$\partial''$ is the remaining part of $\partial$. We have to define
this new grading in such a way that the new differential $\partial'$
is well defined and the corresponding homology groups $H({\cal
C}(K),\partial')$ are invariant with respect to Reidemeister moves.
\bigbreak

A Khovanov homology theory for virtual knots has been constructed
in a sequence of papers by Manturov. In \cite{AddGr}, one gives a
certain procedure for further generalization of these invariants,
which deals with so-called {\it dotted gradings}. In working with Khovanov homology we use
{\it enhanced states} of the Kauffman bracket polynomial. These enhanced states are collections of
labelled simple closed curves obtained by smoothing crossings in the diagram. Each curve is labelled with either
the algebra element $X$ or with the number $1$. The elements $X$ and $1$ belong to the algebra
$k[X]/(X^{2})$ where $k = Z[A,A^{-1}].$
In the dotted grading,   the $X$ and the $1$ can acquire a dot in the form $\dot{X}$ and $\dot{1}.$
We explain how this notation works in the discussion below.
 \bigbreak

We assume all circles in Kauffman states of a diagram can be
assigned a mod ${\bf Z}_{2}$ {\em dotting}: every state circle is
either dotted or not (the dotting should be read from the
topology/combinatorics of the diagram), and the new integral grading
of a chain is set to be ${\#}{\dot X}-{\#}{\dot 1}$, i.e. the number
of dotted circles with the element $X$ minus the number of dotted
circles carrying the element $1$. If this
dotting satisfies certain very simple axioms \cite{AddGr}, then the complex is
well defined and its homology is invariant under Reidemeister moves.
\bigbreak

Another way to introduce the gradings for a given Khovanov homology
theory is to take the coefficients like $[S]$ or $ \prod_{c(S)}
K_{n(s)}$ to be new (multi)gradings themselves, but for this we use
${\bf Z}_{2}$-coefficients. Possibly, this ${\bf
Z}_{2}$-reduction can be avoided if we use twisted
coefficients similar to those from \cite{Izv}, but this has not been done so far. \bigbreak

We note that this paper makes use of enhanced states of the bracket polynomial for discussing
Khovanov homology. This approach was introduced in \cite{Viro1,Viro2}. The first categorification of link invariants in thickened surfaces, thus also of the Kauffman bracket of virtual
links occurs in \cite{APS}. Finally, two recent papers by\cite{Caprau,CMW} can also be viewed as categorifying
the arrow polynomial. although that was not the principle aim of these works. A sequel to this paper
will discuss these relationships.
\bigbreak

\subsection{Acknowledgements}

The last two authors (L.H.K. and V.O.M.) express their gratitude to
the Mathematisches Forschungsinstitut-Oberwolfach, where this paper
was finished, for a nice creative scientific atmosphere.

\section{The Arrow Polynomial $\langle K \rangle_{A}$}
In this section we describe the arrow polynomial invariant
\cite{ExtBr,DK}. One way to see the definition of the arrow
polynomial is to begin with the extended bracket invariant
\cite{ExtBr} and simplify it. The extended invariant is a sum of
graphs (taken up to virtual equivalence in the plane) weighted by
polynomials. In the extended bracket one uses an oriented expansion
so that the smoothings consist of oriented smoothings and
disoriented smoothings. At a disoriented smoothing one sees two
cusps with orientation arrows going into the cusp point in one cusp
and out of the cusp point for the other cusp. Rules for reducing the
states of the extended bracket keep the cusps paired whenever
possible. If we release the cusp pairings at the disoriented smoothings,
we get simpler graphs. These are
composed of disjoint collections of circle graphs that are labelled
 with the orientation markers and
left-right distinctions that occur in the state expansion.
 The basic
conventions for this  simplification are shown in Figure \ref{Figure
2}. In that figure we illustrate how the disoriented smoothing is a
local disjoint union of two vertices (the cusps). Each cusp is
denoted by an angle with arrows either both entering the cusp or
both leaving the cusp. Furthermore, the angle locally divides the
plane into two parts: One part is the span of an acute angle (of
size less than $\pi$); the other part is the span of an obtuse
angle. We refer to the span of the acute angle as the {\it inside}
of the cusp. In Figure \ref{Figure 2}, we have labelled the insides of
the cusps with the symbol $\sharp .$  \bigbreak

\begin{figure}
     \begin{center}
     \begin{tabular}{c}
     \includegraphics[width=6cm]{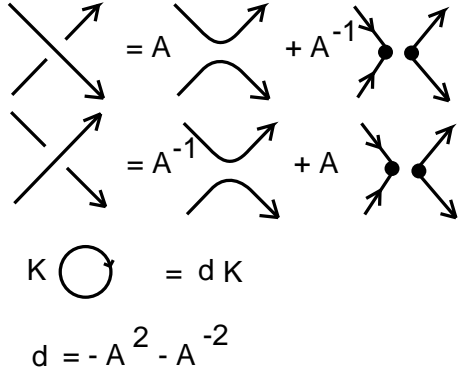}
     \end{tabular}
     \caption{\bf Oriented Bracket Expansion.}
     \label{Figure 1}
\end{center}
\end{figure}

\begin{figure}
     \begin{center}
     \begin{tabular}{c}
     \includegraphics[width=6cm]{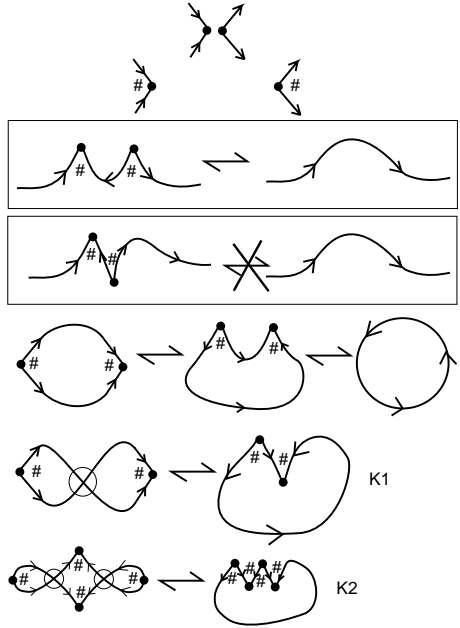}
     \end{tabular}
     \caption{\bf Reduction Relation for the Arrow Polynomial.}
     \label{Figure 2}
\end{center}
\end{figure}

\begin{figure}
     \begin{center}
     \begin{tabular}{c}
     \includegraphics[width=6cm]{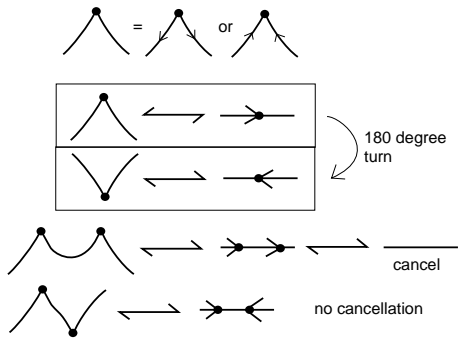}
     \end{tabular}
     \caption{\bf Arrow Convention.}
     \label{Figure 3}
\end{center}
\end{figure}

Figure \ref{Figure 1} illustrates the basic oriented bracket
expansion formula. Figure \ref{Figure 2} illustates the reduction
rule for the {\it arrow polynomial.}  While we have indicated (above) the relationship of the
arrow polynomial with the extended bracket polynomial, the reduction rule for
the arrow polynomial is completely described by Figure \ref{Figure 2}.
We shall denote the arrow polynomial by the notation $\langle K \rangle_{A},$ for a virtual knot or link diagram $K.$ The
reduction rule allows the cancellation of two adjacent cusps when
they have {\it insides on the same side} of the segment that
connects them. When the insides of the cusps are on opposite
sides of the connecting segment,
then no cancellation is allowed. All graphs are taken up to virtual equivalence.
Figure \ref{Figure 2} illustrates the simplification of two circle graphs. In one
case the graph reduces to a circle with no vertices. In the other
case there is no further cancellation, but the graph is equivalent
to one without a virtual crossing. The state expansion for $\langle K \rangle_{A}$ is
exactly as shown in Figure \ref{Figure 1}, but we use the reduction
rule of Figure \ref{Figure 2} so that each state is a disjoint union
of reduced circle graphs. Since such graphs are planar, each is
equivalent to an embedded graph (no virtual crossings) and the
reduced forms of such graphs have $2n$ cusps that  alternate in
type around the circle so that $n$ are pointing inward and $n$ are
pointing outward. The circle with no cusps is evaluated as $d =
-A^2 - A^{-2}$ as is usual for these expansions and the circle is
removed from the graphical expansion. Let $K_{n}$ denote the circle
graph with $2n$ alternating vertex types as shown in Figure
\ref{Figure 2} for $n=1$ and $n=2.$ By our conventions for the
extended bracket polynomial, each circle graph contributes $d = -A^2
- A^{-2}$ to the state sum and the graphs $K_{n}$ (with $n \geq 1$)
remain in the graphical expansion. For the arrow polynomial $\langle K \rangle_{A}$
we can regard each $K_{n}$ as an extra variable in the polynomial.
Thus a product of the $K_{n}$'s denotes a state that is a disjoint
union of copies of these circle graphs with multiplicities. By
evaluating each circle graph as $d = -A^2 - A^{-2}$ we guarantee
that the resulting polynomial will reduce to the original bracket
polynomial when each of the new variables $K_{n}$ is set equal to
unity. Note that we continue to use the caveat that an isolated
circle or circle graph (i.e. a state consisting in a single circle
or single circle graph) is assigned a loop value of unity in the
state sum. This assures that $\langle K \rangle_{A}$ is normalized so that the unknot
receives the value one. \bigbreak

\noindent Formally, we have the following state summation for the arrow polynomial $$ \langle K \rangle_{A} = \sum_{S} \langle K|S\rangle d^{||S||-1} P[S]$$
where $S$ runs over the oriented bracket states of the diagram, $\langle K|S \rangle$ is the usual product of vertex weights as in the
standard bracket polynomial, $||S||$ is the number of circle graphs in the state $S$, and $P[S]$ is a product of the variables $K_{n}$ associated
with the non-trivial circle graphs in the state $S.$ Note that each circle graph (trivial or not) contributes to the power of $d$ in the state summation,
but only non-trivial circle graphs contribute to $P[S].$ The regular isotopy invariance of $\langle K \rangle_{A}$ follows from an analysis of the
behaviour of this state summation under the Reidemeister moves.
\bigbreak

\begin{thm}
With the above conventions, the arrow polynomial
$\langle K \rangle_{A}$ is a polynomial in $A, A^{-1}$ and the graphical variables
$K_{n}$ (of which finitely many will appear for any given virtual
knot or link). $\langle K \rangle_{A}$ is a regular isotopy invariant of virtual
knots and links. The normalized version
$$W[K] = (-A^{3})^{-wr(K)} \langle K \rangle_{A}$$ is an invariant of virtual isotopy.
Here $wr(K)$ denotes the {\it writhe} of the diagram $K$; this is the sum of the signs of
all the classical crossings in the diagram.
If we set $A = 1$ and $d = -A^2 - A^{-2} = -2$, then the resulting specialization
$$F[K] = \langle K \rangle_{A}(A = 1)$$
is an invariant of flat virtual knots and links.
\end{thm}

\bigbreak

 \begin{figure}
     \begin{center}
     \begin{tabular}{c}
     \includegraphics[width=4cm]{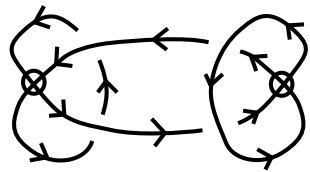}
     \end{tabular}
     \caption{\bf Kishino Diagram.}
     \label{Figure 4}
\end{center}
\end{figure}

\noindent {\bf Example.} Figure \ref{Figure 4} illustrates the
Kishino diagram. With $d= -A^2 - A^{-2}$ $$\langle K \rangle_{A} = 1 + A^4 +
A^{-4} -d^2 K_{1}^{2} + 2 K_{2}.$$  Thus the simple extended bracket
shows that the Kishino is non-trivial and non-classical. In fact,
note that
$$F[K] = 3 + 2 K_{2} - 4 K_{1}^{2}.$$ Thus the invariant $F[K]$ of
flat virtual diagrams proves that the flat Kishino diagram is
non-trivial. This example shows the power of the arrow polynomial.
See \cite{ExtBr,DK} for the details of this calculation. \bigbreak

\section{Khovanov homology for virtual knots}

In this section, we describe Khovanov homology for virtual knots
along the lines of \cite{Kh,BN,Izv}.
\bigbreak

The bracket polynomial \cite{KaB} is usually described by the expansion
\begin{equation}
\langle \raisebox{-0.25\height}{\includegraphics[width=0.5cm]{skcrossr.eps}} \rangle=A \langle \skcrh \rangle + A^{-1}\langle
\skcrv \rangle \label{kabra}
\end{equation}
Letting $c(K)$ denote the number of crossings in the diagram $K,$ if we replace $\langle K \rangle$ by
$A^{-c(K)} \langle K \rangle,$ and then replace $A^{2}$ by $-q^{-1},$ the bracket will be rewritten in the
following form:
\begin{equation}
\langle \raisebox{-0.25\height}{\includegraphics[width=0.5cm]{skcrossr.eps}} \rangle=\langle \skcrh \rangle -q \langle
\skcrv \rangle \label{kabr}
\end{equation}
with $\langle \bigcirc \rangle=(q+q^{-1})$.
In this form of the bracket state sum, the grading
of the Khovanov homology (which is described below) appears naturally. We shall
continue to refer to the smoothings labelled $q$ (or $A^{-1}$ in the
original bracket formulation) as {\it $B$-smoothings}. We should
further note that we use the well-known convention of {\it enhanced
states} where an enhanced state has a label of $1$ or $X$ on each of
its component loops. We then regard the value of the loop $(q + q^{-1})$ as
the sum of the value of two circles: a circle labelled with a $1$ (the value is
$q$) and a circle labelled with an $X$ (the value
is $q^{-1}).$ \bigbreak

To see how the Khovanov grading arises, consider the form of the expansion of this version of the
bracket polynomial in enhanced states. We have the formula as a sum over enhanced states $s:$
$$\langle K \rangle = \sum_{s} (-1)^{n_{B}(s)} q^{j(s)}$$
where $n_{B}(s)$ is the number of $B$-type smoothings in $s$, $\lambda(s)$ is the number of loops in $s$ labelled $1$ minus the number of loops
labelled $X,$ and $j(s) = n_{B}(s) + \lambda(s)$.
This can be rewritten in the following form:
$$\langle K \rangle = \sum_{i, \,j} (-1)^{i} q^{j} \left[ \sum_{s: {n_{B}(s)=i, \, j(s) = j}} 1 \right] =  \sum_{i \,,j} (-1)^{i} q^{j} dim({\cal C}^{ij}).$$
\bigbreak

In the Khovanov homology, the states with $n_{B}(s)=i$ and $ j(s) = j$ form the basis for a module ${\cal C}^{ij}$ over the ground ring
$k.$ Thus we can write $$dim({\cal C}^{ij}) =  \sum_{s: {n_{B}(s)=i, \, j(s) = j}} 1.$$  The bigraded complex composed of the ${\cal C}^{ij}$ has a
differential
$d:{\cal C}^{ij} \longrightarrow {\cal C}^{i+1 \, j}.$ That is, the differential increases the {\it homological grading} $i$ by $1$ and preserves the
{\it quantum grading} $j.$ Below, we will remind the reader of the formula for the differential in the Khovanov complex. Note however that the
existence of a bigraded complex of this type allows us to further write:
$$\langle K \rangle = \sum_{j} q^{j} \sum_{i} (-1)^{i} dim({\cal C}^{ij}) = \sum_{j} q^{j} \chi({\cal C}^{\bullet \,
j}),$$
where $\chi({\cal C}^{\bullet \, j})$ is the Euler characteristic of the subcomplex ${\cal C}^{\bullet \, j}$ for a fixed value of $j.$ Since $j$ is
preserved by the differential, these subcomplexes have their own Euler characteristics and homology. We can write
$$\langle K \rangle = \sum_{j} q^{j} \chi(H({\cal C}^{\bullet \, j})),$$ where $H({\cal C}^{\bullet \, j})$ denotes the homology of this complex.
Thus our last formula expresses the bracket polynomial as a {\it graded Euler characteristic} of a homology theory associated with the enhanced states
of the bracket state summation. This is the categorification of the bracket polynomial. Khovanov proves that this homology theory is an invariant
of knots and links, creating a new and stronger invariant than the original Jones polynomial.
\bigbreak

We explain the differential in this complex for mod-$2$
coefficients and leave it to the reader to see the references for
the rest. The differential is defined via the algebra ${\cal A}=
k[X]/(x^{2})$ so that $X^2 = 0$ with coproduct $\Delta: {\cal A}
\longrightarrow {\cal A} \otimes {\cal A}$ defined by $\Delta(X) =
X \otimes X$ and $\Delta(1) = 1 \otimes X + X \otimes 1.$ Partial
differentials (which are defined on an enhanced state with a
chosen site, whereas the differential is a sum of these mappings)
are defined on each enhanced state $s$ and a site $\kappa$ of type
$A$ in that  state. We consider states obtained from the given
state by  smoothing the given site $\kappa$. The result of
smoothing $\kappa$ is to produce a new state $s'$ with one more
site of type $B$ than $s.$ Forming $s'$ from $s$ we either
amalgamate two loops to a single loop at $\kappa$, or we divide a
loop at $\kappa$ into two distinct loops. In the case of
amalgamation, the new state $s$ acquires the label on the
amalgamated circle that is the product of the labels on the two
circles that are its ancestors in $s$. That is, $m(1 \otimes X) =
X$ and $ m(X \otimes X) = 0 $. Thus this case of the partial
differential is described by the multiplication in the algebra. If
one circle becomes two circles, then we apply the coproduct. Thus
if the circle is labelled $X$, then the resultant two circles are
each labelled $X$ corresponding to $\Delta(X) = X \otimes X$. If
the orginal circle is labelled $1$ then we take the partial
boundary to be a sum of two enhanced states with  labels $1$ and
$X$ in one case, and labels $X$ and $1$ in the other case  on the
respective circles. This corresponds to $\Delta(1) = 1 \otimes X +
X \otimes 1.$ Modulo two, the differential of an enhanced state is
the sum, over all sites of type $A$ in the state, of the partial
differential at these sites. It is not hard to verify directly
that the square of the differential mapping is zero and that it
behaves as advertised, keeping $j(s)$ constant. There is more to
say about the nature of this construction with respect to
Frobenius algebras and tangle cobordisms. See \cite{Kh,BN,BN2}
\bigbreak

Here we consider bigraded complexes ${\cal C}^{ij}$ with {\em
height} (homological grading) $i$ and {\em quantum grading} $j.$  In the
unnormalized Khovanov complex $[[K]]$ the index
$i$ is the number of $B$-smoothings of the bracket, and for
every enhanced state, the index $j$ is equal to the number of components labelled $1$ minus the number
of components labelled $X$ plus the number of $B$-smoothings. The normalized complex
differs from $[[K]]$ by an overall shift of both gradings; the
differential preserves the quantum grading and increases the height
by $1$. The height and grading shift operations are defined as
$({\cal C}[k]\{l\})^{ij}={\cal C}[i-k]\{j-l\}$.\label{shifts}
\bigbreak

This form is used as the starting point for the Khovanov homology.
We now describe the formalism in a bit more detail in order to give the structure of the differential for
Khovanov homology of virtual knots and links.
For a diagram $K$ of a virtual knot, we consider the {\em state cube} defined as follows:
Enumerate all $n$ classical crossings of $K$ in arbitrary way and consider all Kauffman states
(states as collections of loops without specific enhancement labels) as vertices of the discrete cube $\{0,1\}^{n}.$
Each coordinate corresponds to a way of smoothing and is equal to
$0$ (the $A$-smoothing) or $1$ (the $B$-smoothing). Thus, each vertex of the cube defines a set of circles
(say, $p$ circles), and this set of circles defines a certain vector space (module) of dimension $2^{p}.$ The module for a single
circle is generated by $1$ and $X.$ The spaces together form the total chain space of the unnormalized Khovanov complex $[[K]]$
and its normalized version ${\cal C}[[K]]$.
We omit the normalisation, which is standard, and refer the reader to \cite{Kh,BN,Izv}.
\bigbreak

We regard the loop factors for the unenhanced bracket, $(q+q^{-1}),$ as graded dimensions of
the module $V=Span(\{1,X\}), \mbox{deg}\; 1=1, \mbox{deg}\; X=-1$ over some ring $k$,
and the height $i(s)$ plays the role of homological dimension.
Define the chain space $[[K]]_{i}$ of homological
dimension $i$ to be the direct sum over all vertices of height $i$ (defined
as above) of
$V^{\gamma(s)}\{i\}$ (here $\{\cdot\}$ is the quantum grading
shift and $\gamma(s)$ is the number of loops in the state $s$). Then the alternating sum of graded dimensions of
$[[K]]_{i}$, is precisely equal to the (modified) Kauffman bracket, as we have described above.
\bigbreak

Thus, if one defines a differential on $[[K]]$ that preserves the
grading and increases the homological dimension by $1$, the Euler
characteristic of that complex will be precisely the
bracket.
\bigbreak

We now consider a generalization of the Khovanov homology to virtual knots. When we pass from one state of the state cube to a neighboring state (which
differs precisely at one coordinate), we get a resmoothing of the set of circles.
We refer to that as a {\em bifurcation} of the state cube. Such a bifurcation can either merge two circles into one ($2\to 1$-bifurcation)
or split one circle
into two ($1\to 2$-bifurcation), or (in the case of virtual knots and links) transform one circle into one ($1\to 1$-bifurcation).
These bifurcations encode the information about differentials in the complex as follows.
\bigbreak

We have defined the {\em state cube} consisting of state loops and
carrying no information how these loops interact. For
Khovanov homology, we deal with the same cube, remembering the
information about the loop bifurcation. Later on, we refer to it
as a {\em bifurcation cube}.
\bigbreak

The chain spaces of the complex are well defined. However, the
problem of finding a differential $\partial$ in the general case of
virtual knots, is not easy. See Figure 7 for a key example that we shall discuss.
To define the differential, we have
to pay attention to the different isomorphism classes of the chain space
identified by using local bases (see below).
\bigbreak

The differential acts on the chain space as follows: it takes a
chain (regard an enhanced state as an elementary chain) corresponding to a certain vertex of the bifurcation cube to
some chains corresponding to all adjacent vertices with greater
homological degree. That is, the differential is a sum of {\em
partial differentials}, each partial differential acts along an edge
of the cube. Every partial differential corresponds to some
direction and is associated with some classical crossing of the
diagram. The total differential is the sum of these partial differentials, and
so formally looks like $$\partial = \sum_{a} \partial_{a}$$ where the summation is over
all edges of the cube. In discussing differentials we shall often refer to a partial differential
without indicating its subscript.
\bigbreak

Selecting an un-enhanced Kauffman state $S$ (consisting of loops
with cusps), we choose an arbitrary order for the circles in S. and
then orient each circle in $S$. Letting $ \gamma (S) = ||S|| $ be
the number of loops in $S$, associate the module $ \Lambda^{ ||S||}
(V) $ to $S$ where this denotes the $ ||S||^{th} $ exterior power of
$V$ -- the order of the factors in the exterior power depends on the
choice of the ordering that was chosen. Having made this choice (of
ordering and orientation), if $s$ is an enhancement of $S$ then
label all loops in the state $s$ with either $ +X $ or $+1$
according to the enhancement. This oriented, ordered, and labelled
state forms a generating chain in the complex. If the orientation of
a loop in $S$ is reversed then the label for $ X $ becomes $-X$ but
the label for $1$ does not change. Otherwise, signs change according
to the structure of the exterior algebra.

Then for a state with $l$ circles, we get a vector space (module) of
dimension $2^{l}$. All these chains have homological dimension $i = n_{B}$.
We set the quantum grading $j$ of these chains equal to  $i$ plus the number of
circles marked by $\pm 1$ minus the number of circles marked by $\pm
X$.
\bigbreak

\begin{figure}
\centering\includegraphics[width=300pt]{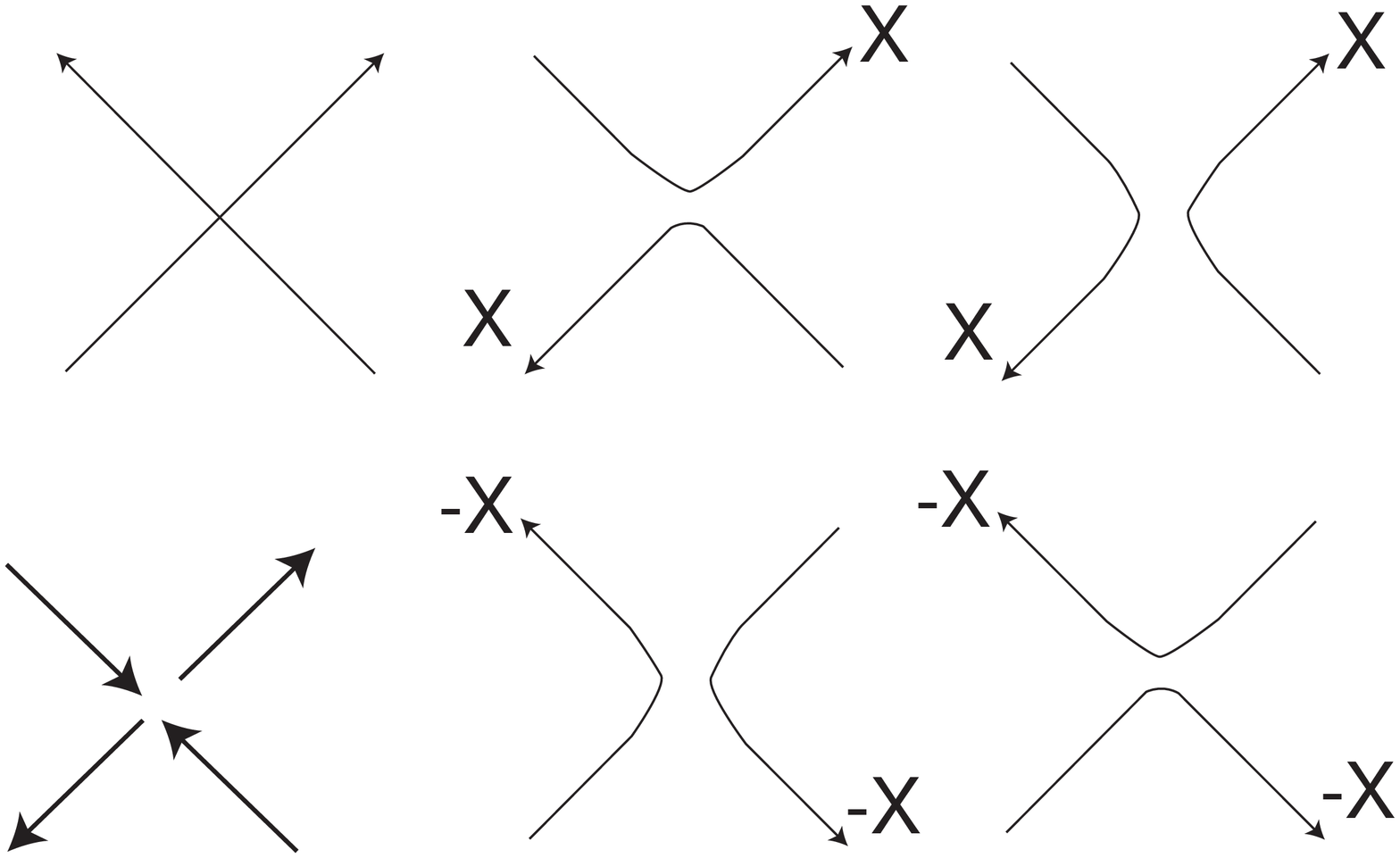} \caption{Setting
the local basis for a crossing} \label{figj}
\end{figure}

Let us now define the partial differentials of our complex. First,
we think of each classical crossing so that its edges are oriented
upwards, as in Figure \ref{figj}, upper left picture.
\bigbreak

Choose a certain state of a virtual link diagram $L\subset {\cal
M}$. Choose a classical crossing $U$ of $L$. We say that in a state
$s$ that a state circle $\gamma$ is incident to a classical crossing $U$
if at least one of the two local parts of smoothed crossing $U$
belongs to $\gamma$. Consider all circles
 $\gamma$ incident to  $U$. Fix some orientation of these circles according to
the orientation of the edge emanating in the upward-right direction
and opposite to the orientation of the edge coming from the bottom
left, see Figure \ref{figj}. Such an orientation is well defined
except for the case when resmoothing one edge takes one circle to one circle. In
such a situation, we shall
not define the local basis  $\{1,X\}$, and we set the partial
differential corresponding to that edge to be zero.
\bigbreak

In the other situations, the edge of the cube corresponding to the
partial differential either increases or decreases the number of
circles. This means that at the corresponding crossing the local
bifurcation either takes two circles into one or takes one circle
into two. If we deal with two circles incident to a crossing from
opposite signs, we order them in such a way that the upper (resp.,
left) one is the first one; the lower (resp., right) one is the
second; here the notions ``left, right, upper, lower'' are chosen
according to the rule for identifying the crossing neighbourhood
with Figure \ref{figj}. Furthermore, for defining the partial
differentials of types $m$ and $\Delta$ (which correspond to
decreasing/increasing the number of circles by one) we assume that
the circles we deal with are in the initial positions specified in our
ordered tensor product; this can always be achieved by a preliminary
permutation, which, possibly leads to a sign change. Now, let us
define the partial differential locally according to the prescribed
choice of generators at crossings and the prescribed ordering.
\bigbreak

Now, we describe the partial differentials $\partial$ from
\cite{Izv} without new gradings. If we set $\Delta (1)=1_{1}\wedge
X_{2}+X_{1}\wedge 1_{2}; \Delta(X)=X_{1}\wedge X_{2}$ and
$m(1_{1}\wedge 1_{2})=1;m(X_{1}\wedge 1_{2})=m(1_{1}\wedge X_{2})=X;
m(X_{1}\wedge X_{2})=0$, define the partial differential $\partial$
according to the rule $\partial (\alpha\wedge \beta)=m(\alpha)\wedge
\beta$ (in the case we deal with a $2\to 1$-bifurcation, where
$\alpha$ denotes the first two circles $\alpha$) or
$\partial (\alpha\wedge \beta)=\Delta(\alpha)\wedge \beta$ (when one
circle marked by $\alpha$ bifurcates to two ones); here by $\beta$
we mean an ordered set of oriented circles, not incident to the
given crossings; the marks on these circles $\pm 1$ and $\pm X$ are
given.
\bigbreak

\begin{thm}\cite{Izv}
Let $K$ be a virtual knot or link. Then
$[[K]]$ is a well-defined complex with respect to $\partial.$ After
a small grading shift and a height shift, the homology of $[[K]]$ is invariant
under the generalised Reidemeister moves for virtual knots and links.
\end{thm}

\section{Grading Considerations for the Arrow Polynomial $\langle K \rangle_{A}$}
\begin{figure}
     \centering\includegraphics[width=6cm]{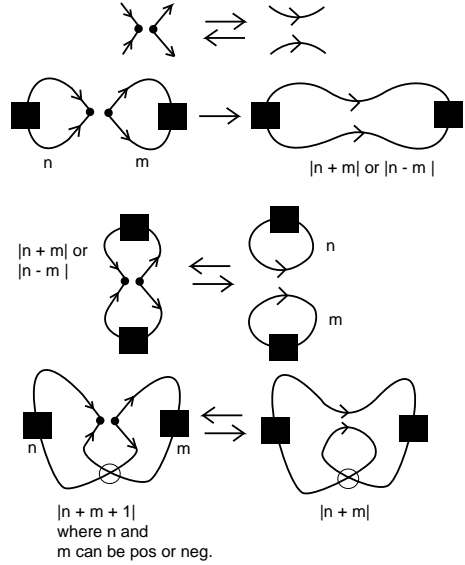}
     \caption{\bf Arrow Numbers for Interacting Loops}
\label{Figure6}
\end{figure}

In order to consider gradings for Khovanov homology in relation to
the structure of the arrow polynomial $ \langle K\rangle_{A}$ we have to examine how the
arrow number of state loops change under a replacement of an
$A$-smoothing by a $B$-smoothing. Such replacement, when we use
oriented diagrams involves the replacement of a cusp pair by an
oriented smoothing or vice versa. Furthermore, we may be combining
or splitting two loops. Refer to Figure \ref{Figure6} for a
depiction of the different cases. This figure shows the three basic
cases. \bigbreak

In the first case we have two loops $C_{1}$ and $C_{2}$ sharing a disoriented site and the
smoothing is  a single loop $C$ where the paired cusps of the disoriented site disappear. In this case if $n(C_{1}) = n$ and $n(C_{2}) = m$, then
 $n(C) = |n-m|.$
\bigbreak

In the second case, we have a single loop $C$ with a disoriented site and a pair of cusps, and on smoothing this site we obtain two loops $C_{1}$ and
$C_{2}$ whose arrow numbers are $n(C_{1}) = n$ and $n(C_{2}) = m.$ The following arrow numbers for $C$ are then possible
$|n(C)| = |n-m|$ or $|n+m|.$
\bigbreak

In the third case, we have a single loop $C$ with a disoriented site and a pair of cusps, and on smoothing this site we obtain a single loop $C'$.
Assuming that $n(C') = |n + m|$ as shown in the figure, we have $|n(C)| = |n+m + 1|$ where $n$ and $m$ can be positive or negative.
\bigbreak

These are all the ways that loops can interact and change their respective arrow numbers. In the next section, we will apply these results to the
grading in Khovanov homology.
\bigbreak

\section{Dotted gradings and the dotted categorification}

First, we introduce a concept of {\em dotting axiomatics} as developed in \cite{AddGr}.
The purpose of this dotting axiomatics is to give general conditions under which extra decorations on the states can be used
to create new gradings and hence new versions of Khovanov homology. We will apply these axiomatics to the arrow numbers
on the state loops of the arrow polynomial.
\bigbreak

For the axiomatics, assume we have some class of
objects with Reidemeister moves, Kauffman bracket and the Khovanov
homology (in the usual setup or in the setup of \cite{Kho}).
Assume that there is a method, which for every diagram and every
state of it associates dots to some of the circles in the bracket states in such a way
that the following conditions hold:

\begin{enumerate}

\item The dotting of circles is additive with respect to $2\to
1$-bifurcations and $1\to 2$-bifurcations mod $2$. This additivity means that
when we merge two circles (split one circle into two), the number of dots on the circles being operated
on is preserved modulo ${\bf Z}_{2}$.

This means that the parity of the number of dots on the circles operated on is preserved
whenever we merge two circles or split one circle into two.

If the dotting is not preserved under a $ 1 \to 1$ bifurcation, then this bifurcation is
taken to be the zero map.

\item Similar curves for corresponding  smoothings of the RHS and the LHS
of any Reidemeister move have the same dotting.

\item Small circles appearing for the first, the second, and the
third Reidemeister moves are not dotted.

\end{enumerate}
\bigbreak

Let us call the conditions above {\em the dotting conditions}.
With such a structure in hand, one defines a {\it new grading $g(s)$} for states $s$ by taking the difference between
the number of dotted $X$'s and the number of dotted $1$'s in the state.
$$g(s) = \sharp ({\dot X}) - \sharp ({\dot 1})$$
We shall use this grading in the
constructions that follow.
\bigbreak

\begin{thm}
Assume there is a theory using the Khovanov complex $([[K]],\partial)$
such that the Kauffman states can be dotted so that the dotting
conditions hold. Take $[[K]]_{g}$ to be the space $[[K]]$ endowed with new grading
as above.
\bigbreak

Define $\partial'$ to be the composition of $\partial$ with the new grading projection
and set $\partial''=\partial-\partial'$.
\bigbreak

Then the homology of $[[K]]_{g}$ (with respect to $\partial'$) is invariant (up to a degree shift
and a height shift).
\bigbreak

For any operator $\lambda$ on the ground ring, the complex
$[[K]]_{g}$ is well defined with respect to the differential
$\partial'+\lambda\partial''$, and the corresponding homology is
invariant (up to well-known shifts).
\bigbreak

Moreover, if we have several forms of dotting $g_{1}, g_{2}, \dots, g_{k}$
occuring together on the same Khovanov complex
so that for each of them the dotting condition holds, then the
complex $K_{g_{1},\dots,g_{k}}$ with differential
$\partial_{g_{1},\dots, g_{k}}$ defined to be the projection of
$\partial$ to the subspace preserving all the gradings, is
invariant.
\bigbreak

\label{mainthm}
\end{thm}

The theorem above allows one to `raise' some additional information modulo ${\bf Z}_{2}$ to
the level of gradings. Our aim is to categorify the arrow polynomial, that is, to add new gradings corresponding
to the arrow count: for every state we have a set of circles labelled by a set of non-zero integers, and this set of integers should
be represented in the complex as a grading. Theorem \ref{mainthm} shows that it is possible to do that when we consider the information of the arrow count
only modulo ${\bf Z}_{2}$: the conditions of additivity and similarity under Reidemeister moves for arrow count were checked in the previous section
of this paper.
\bigbreak

In order to use the integral information about the arrow count, we have to undertake a generalization of the construction of
theorem \ref{mainthm}. We shall do this in the next section. This section of the paper is devoted to describing a first-order
categorification of the arrow polynomial.
\bigbreak

The main idea behind the proof of Theorem \ref{mainthm} is as follows.
Additivity of the grading can be verified and checked
on a bifurcation cube. First of all, it follows from a straightforward
check that $\partial''$
always increases the dotted grading (this is proved in \cite{Izv} but can be taken here as an exercise for the reader).
Then, the complex is well defined because $(\partial')^{2}$ is
nothing but a composition of $(\partial)^{2}$ with a
``grading-preserving projection''. This is guaranteed because
$\partial''$ strictly increases the new grading. Note the the mod-2 preservation of the dotting is what makes this grading increase
of $\partial''$ work. thus Theorem \ref{mainthm} depends ultimately on that parity presevation of the dotted grading.
\bigbreak

The main idea of the invariance under Reidemeister moves is similar
to the usual Khovanov idea, see for example \cite{BN}: we have to check
that the multiplication $m$ remains surjective after reducing
$\partial$ to $\partial'$ and $\Delta$ remains injective. The latter
follows from the fact that ``small circles are not dotted''.
\bigbreak

Now, one can easily check that the conditions of the theorem hold if
we set the dotting as follows: the curve is dotted if it is marked as $K_{j}$ with $j$ odd,
and it is not dotted if it is marked as $K_{i}$ with $i$ even.
\bigbreak

Now, one checks that

\begin{figure}
\centering\includegraphics[width=300pt]{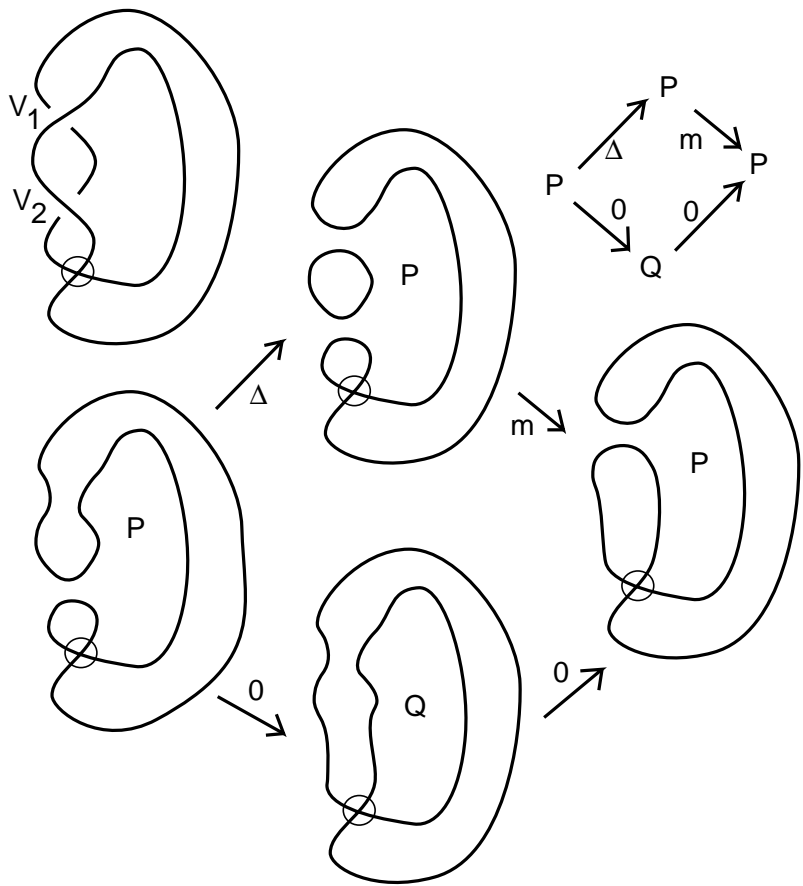} \caption{A
face for the case of the extended bracket of the cube}
\label{norient}
\end{figure}

\begin{enumerate}

\item The dotting is ${\bf Z}_{2}$-additive with respect to
resmoothing (performing $1\to 2$ or $2\to 1$ bifurcation).
\bigbreak

This follows from Figure \ref{Figure6} upper part: we see that when
merging two circles with arrow count $m$ and $n$, we get $\pm m\pm
n$ and when splitting a circle with arrow number $k$, we get two
circles with arrow numbers $l$ and $\pm k\pm l$ which results in
${\bf Z}_{2}$-additivity under $2\to 1$ and $1\to 2$-bifurcations.
\bigbreak

On the other hand, if partial differentials for all $1\to 1$
bifurcations are set to be zero, it can be checked that all faces
having at least $1\to 1$-bifurcation are anticommutative because
$0=0$. The only non-trivial example is shown in Figure \ref{norient},
and the corresponding calculation is performed in \cite{Izv}.
\bigbreak

\item The small circles coming from Reidemeister moves are not
dotted.
Indeed, for the 1st Reidemeister move we have no cusps at all,
and for the second move and for the third move we have two cusps of
opposite signs.
\bigbreak

\item For any Reidemeister move, the corresponding state diagrams in the LHS and
RHS have the same dotting. Locally, there is no grading change for the
Reidemeister moves when we use arrow counts.
Again, this follows from the {\it invariance} under Reidemeister
moves: two pictures would not get cancelled if they had different
coefficients coming from cusps; this means they have the same
dotting.

\end{enumerate}

\bigbreak

\section{${\bf Z}_{2}$-categorification with general gradings}

\subsection{General setup}

The aim of this section is to prove a general theorem on
categorification that fits the arrow polynomial. This is an
extension of the dotted grading construction, which works, however,
only with ${\bf Z}_{2}$-coefficients for the homology. Later, we shall discuss
whether this construction can be extended to the case of integral
coefficients. For instance, we can extend this construction to the
case of integral coefficients if the odd Khovanov
homology theory \cite{ORS} can be defined for this class of knots.
\bigbreak

Briefly, we want to start with a Khovanov homology (usual over ${\bf
Z}_{2}$ or the one using twisted coefficients) and make some partial
differentials equal to zero.
\bigbreak

As the initial data for this theorem, we require that
we have a well-defined bracket, and we assume
that in each state of the diagram, each circle is given a
non-negative integer. For the dotted conditions, we require
that

\begin{enumerate}

\item The numbers are ``plus-minus additive'' with respect to $2\to
1$-bifurcations and $1\to 2$-bifurcations, that is, if a resmoothing
of two circles labelled by non-negative integers $p$ and $q$ leads to
one circle, the label of this circle will be $|p+q|$ or $|p-q|$.

\item Similar curves for corresponding smoothings of the RHS and the LHS
of any Reidemeister move have the same numbers.

\item Small circles appearing for the first, the second, and the
third Reidemeister moves are labelled by zeroes.

\end{enumerate}

We call these conditions {\em integer labelling conditions}.
\bigbreak

After this, our strategy will be as follows:
If we attempt to make the integral arrow count the new grading
and take the part of the differential preserving this, to be the new differential, we shall
see that the square of this new differential will not be zero. Consider the situation when
a $2$-face of the bifurcation cube has arrow counts $P$ in the left corner (both smoothings
zero), $P$ in the upper corner (both smoothings one), $P$ in one right corner and $Q$
in the remaining corner. See Figure 7 for an example. Then one composition of the two differentials (going through $P$) survives,
while the other one (going through $Q$) becomes zero. That is why the square of the new proposed differential,
detecting the arrow count, is non-zero. On the other hand, all the information about the arrow count has to
be included in order to get a faithful categorification of the arrow polynomial (that is, having a
chain space with gradings one can restore the arrow count, and having the homology one can restore the
arrow polynomial). In order to solve this problem, we are going to introduce two new sorts of gradings, one
of which will correct the other, and make the differential well-defined.
\bigbreak

We take the usual
Khovanov differential $\partial$ and form two new series of gradings
(called {\em multiple gradings} and {\em vector gradings}).
After that, for each basic chain of the complex we have a whole
collection of gradings, and we define the new differential
$\partial'$ to be the composition of $\partial$ with the projection
to the subspace where all gradings are preserved by $\partial$, having the same gradings
(all multiple and vector gradings) as in the preimage. That
is, we let $S = \{ x| gr(x) = gr(\partial x) \}$ and define
$\partial' = \partial |_{S}.$
\bigbreak

Now, we introduce {\em multiple gradings} as follows. A multiple
grading is a set of strictly positive integers that is associated
with a Kauffman state of the diagram. That is,the state is not yet labelled with $X$ and $1$;
a {\it basic chain} in the state is such a labelling.
With each state, we shall
associate exactly one multiple grading for each basic chain in this
state, independently from the particular choice of $X$ and $1$ on
circles. {\it This multiple grading is just the set of all
non-zero arrow counts on circles of the state.}
\bigbreak

The {\it vector grading} is an infinite ordered collection (list) of
integers (first, second, third, etc.) each of which might be either
positive or negative or zero. The vector grading
depends on the particular choice of $1$ and $X$ on all state circles.
But before introducing the vector grading, we introduce
the {\em vector dotting} for state circles (that have the
initial labelling by arrow numbers).  For a circle labelled by $p$
we put no dots at all if $p=0$; otherwise we represent
$p=2^{k-1}{l}$, where $l$ is odd and put exactly one dot of order
$k$ over this circle (we also call it a {\em $k$-th dot}). Thus, for $p=1$ we will have only one primary
dot, for $p=2$ we will have only one secondary dot, for $p=4$ we
will have only one ternary dot and so on. The vector dotting is an infinite vector
of these {\em dot numbers} with one possibly non-zero coordinate for each state circle. Note that the vector dotting
depends only on arrow numbers for the Kauffman state.
\bigbreak

Now we can define the vector grading.
The vector grading of a trivial circle (without dots) is the zero vector
$(0,\dots,0,\dots )$. For a non-trivial circle having one $k$-th
dot, {\it the grading is set to be $+1$ on $k$-th vector position for the
enhanced state carrying $X$ and $-1$ on $k$-th vector position for the
enhanced state carrying $1$;} the other entries of the vector grading for a given enhanced state circle are set to be zero.
\bigbreak

{\it The vector grading of a basic chain (enhanced state)
is defined to be the coordinatewise sum of the vector gradings (these are infinite
vectors) over all circles in the enhanced state.}
Thus, if we have one circle labelled by $2$ with element
$X$ on it and another circle labelled by $1$ with element $1$ on it,
we get the vector grading: $(-1,1,0,\dots,0,\dots)$.
\bigbreak

The chain space of the initial Khovanov complex is split into
subspaces with respect to the multiple grading and vector
grading. We set the differential $\partial'$ to be the composition
of the initial differential $\partial$ with the projection to the
subspace having the same gradings as the preimage.
\bigbreak

\begin{thm}
If a state labelling satisfies the integer labelling conditions, then the
complex ${\cal C}$ is well defined with respect to differential
$\partial'$ (that is, $(\partial')^{2}=0$), and its homology groups
$H({\cal C}, {\partial'})$ are invariant with respect to the Reidemeister
moves.

\end{thm}

First, let us check that the arrow polynomial statisfies the integer
labelling conditions. This follows from Figure \ref{Figure6}.
Now, the second condition ``similar curves generate similar
smoothing'' also follows from a direct calculation, as well as the
third condition about trivial circles coming from Reidemeister
moves. Indeed, for the first Reidemeister move one gets a small loop
without any cusp, for the second Reidemeister move one gets either a
loop without cusps or a loop with two cusps cancelling each other.
The same for the third Reidemeister move: one gets at least two
cusps, which should cancel each other. This proves that the {\em
integer labelling conditions} hold for the arrow count.
\bigbreak

Now, let us prove the main theorem. The proof will consist of the
two parts: the difficult one, where we show that the complex is well
defined (the square of the differential is zero) and the easy one,
where we prove that the homology is invariant under Reidemeister
moves. The second part will be standard and in main features it will
repeat the analogous proof for the usual Khovanov homology.
\bigbreak

{\bf Part 1. Proof that the complex is well defined.}

We first note that we work over ${\bf Z}_{2}$-coefficients. We have
to prove that for every $2$-face of the bifurcation cube, the two
compositions corresponding to faces will coincide. This means that
commutativity and anticommutativity coincide.
\bigbreak

 An {\em atom} is a pair
 $(M,\Gamma)$ of a $2$-manifold $M$ and a graph $\Gamma$ embedded $M$
 together with a colouring of $M\backslash \Gamma$ in a checkerboard
 manner. Here $\Gamma$ is called the {\em frame} of the atom, whence
by {\em genus} (resp., {\em Euler characteristic}) of the atom we
 mean that of the surface $M$.
\bigbreak

With a virtual knot diagram (with every component having at
least one classical crossing) we associate an atom as follows. (Note that the atom
need not be orientable). We take all
classical crossings to be vertices of the frame. The edges of the frame correspond
to branches of the diagram connecting classical crossings (we do not take
into account how they intersect in virtual crossings).
Moreover, the edges of the frame emanating from a vertex are naturally split
into two pairs of {\em opposite} ones: the opposite relation (ordering of edges)
is taken from the
plane diagram. Thus we get two pairs of opposite edges (opposite in the sense that these edges are not adjacent
in the cyclic order of edges about the vertex) and also four {\em angles}
generated by pairs of adjacent (non-opposite) edges.
Now, for the obtained four-valent
graph we attach black and white cells as follows: for every crossing we indicate two pairs
of adjacent edges for ``pasting the black cells'', and the remaining pair of angles are used for
attaching black cells. Cells are attached globally to
conform these local conditions. The ``black angles'' correspond to pairs of edges taken from the $B$-smoothing
of the bracket. This completely defines the way for attaching black and white cells to
get a $2$-manifold starting from the frame.
\bigbreak

This atomic terminology is useful in classifying virtual diagrams in terms of orientability and non-orientability of the corresponding
atoms. An atom has a $1\to 1$-bifurcation if and only if it is non-orientable \cite{Izv}.
In the following we shall need to discuss all atoms that derive from diagrams with two
crossings. The reader can easily enumerate the possible Gauss codes with two symbols and arrive at the possibilities
$(11)(22)$ (two components, four cases depending on the crossings), $(12)(12)$ (a Hopf link configuration with four crossing possibilities),
$(1122)$ (a single unknotted component), $(1212)$ (a non-orientable atom). These cases need to be analyzed and the reader
will find them depicted in Figures 8. See also Figures 9 through 12.
\bigbreak

Each possible $2$-face of the bifurcation cube represents an atom
with $2$ vertices (that is, the face represents all four possibilities for smoothing a pair of crossings
in the original link diagram): for each atom, there are four states
$AA,AB,BA,BB$ and four maps corresponding to partial differentials $AA\to AB,AA\to
BA,AB\to BB,BA\to BB$. Some of them correspond to $1\to
1$-bifurcation which means that the corresponding partial
differential in the usual Khovanov complex is zero. Thus, so is the
partial differential in question (it is a composition of zero map
with a projection). By parity reasons, for a given atom, there may be
$4, 2$ or $0$ partial differentials (in the initial cube) which are
equal to zero.
\bigbreak

If all four differentials are equal to zero, then we get the desired
equality for the composition of the differentials as $0=0$. If we
have $2$ maps of type $1\to 1$ then two options are possible. In one
of them we have one zero map for each of the two compositions, which
leads to $0=0$. We call such atoms  {\em inessential}. In the other
case we have $0$ for the composition of the two $ 1 \to 1 $ maps,
but the other composition of maps must be analyzed. \bigbreak

\begin{figure} \centering\includegraphics[width=300pt]{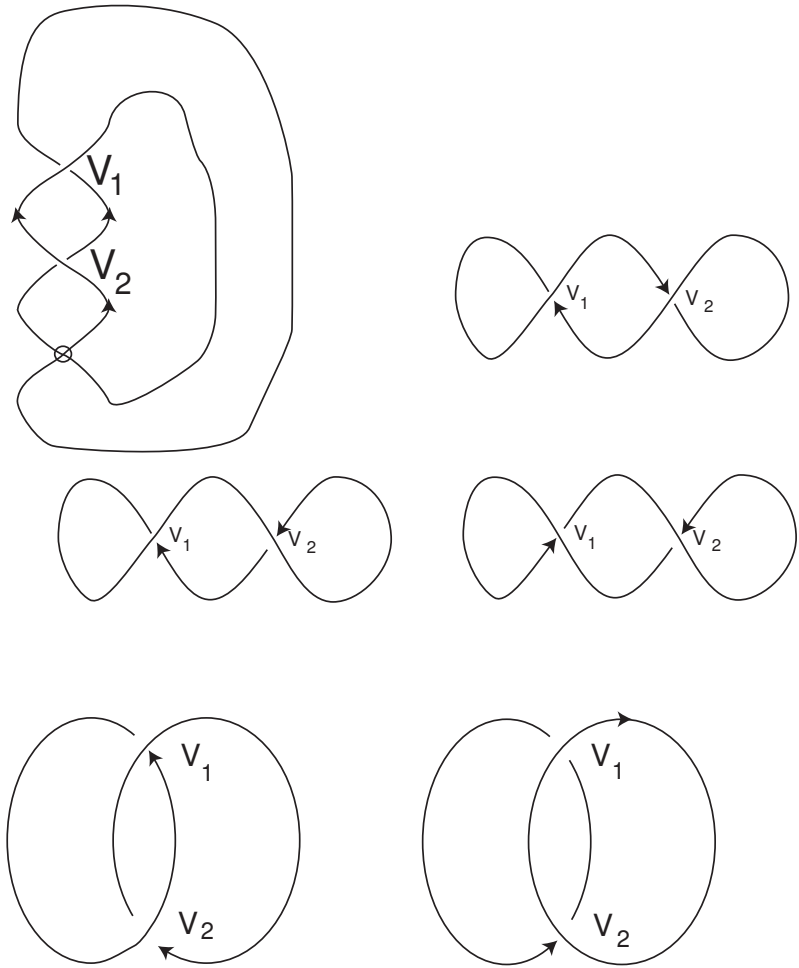}
\caption{Essential atoms with 2 vertices} \label{allatoms}
\end{figure}

Thus we are left with $6$ essential atoms as shown in Figure \ref{allatoms}.
\bigbreak

For each of these atoms the usual Khovanov differential produces a
commutative diagram. Now, multi-gradings and multi-dotted gradings
come into play. We have to show that for each atom $V$ the equality
of partial differentials $q\circ p=s\circ r$ for the usual Khovanov
differentials will hold for the reduced differentials $q'\circ
p'=s'\circ r'$. Here $p,q,r,s$ denote the four partial differentials that occur in the
Khovanov complex at the atom in question.
Some remarks are in order.
\bigbreak

\noindent {\bf Notation.} Let us denote the differential of the Khovanov complex by
$\partial$, and denote its combination with the projection respecting
the multiple grading by $\partial_{multi}$, its combination with the
projection respecting the vector gradings by $\partial_{vect}$
and denote the combination with both projections by $\partial'$. We
are mostly interested in the cases when $\partial=\partial'$ or when
$\partial'=0$ for some particular element of the chain complex.
\bigbreak

We have to list all atoms with two vertices. Some of them are
disconnected in the sense that there is no edge connecting one
vertex to the other.
\bigbreak

For such atoms the (anti)commutativity obviously holds.
\bigbreak

Now, let us list all connected essential atoms. There are exactly
$6$ of them, one non-orientable, $3$ orientable with the frame of
the unlink and $2$ orientable  with the frame of the Hopf link,
see Figure \ref{allatoms}.
\bigbreak

For each atom, the anticommutativity of the virtual Khovanov homology
over ${\bf Z}$ is checked in \cite{Izv}, which leads to the
(anti)commutativity over ${\bf Z}_{2}$. Our goal is to check that
the multigradings and dotted multigradings preserve this
(anti)commutativity.
\bigbreak

For this sake we must consider all possible labellings of the
state circles for atoms. Each labelling gives a number of integers,
for which we take only absolute values and consider only non-zero
ones. This leads to the following multiple gradings $P,Q,R,S$ where
$P$ corresponds to the smoothing of the atom where both crossings
have $A$-type of smoothing, for $S$ both crossings have $B$-type of
smoothing, and for each of $Q,R$ one crossing has $A$-smoothing and
the other one has the $B$-smoothing. See Figure 9.
\bigbreak

\begin{figure}
\centering\includegraphics[width=200pt]{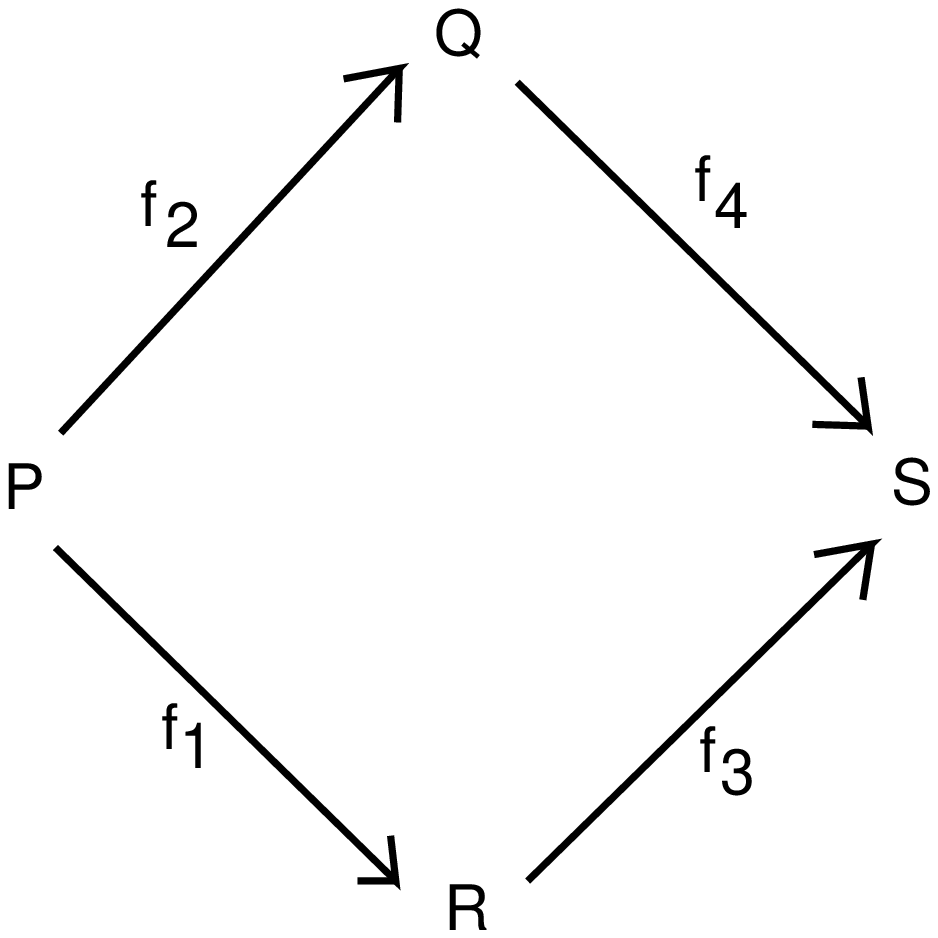}
 \caption{Each atom generates a labelling and two compositions of
maps} \label{atm}
\end{figure}

We must look at the differentials depending on
$P,Q,R,S$. Denote the corresponding partial differentials of $\partial'$
by $f_{1},f_{2},f_{3},f_{4}$, respectively, see Figure \ref{atm}.
\bigbreak

The following lemma holds.

\begin{lm}
If the multiple gradings $P,Q,R,S$ as described above are all equal ($P=Q=R=S$), then for all partial differentials
corresponding to the atom under discussion, we have $\partial_{mult}=\partial$.\label{lmlm1}
\end{lm}
\bigbreak

Let us now look at vector gradings. There is one case when
$\partial_{vect}\neq \partial$ because of the following. Assume we
have a $1\to 2$ or $2\to 1$-bifurcation where all three circles are
dotted: two circles have dotting of order $k$ and one circle has
dotting of higher order $l>k$. This may happen, e.g., in $2\to
1$-bifurcation, when the two circles to merge have arrow label one each
(one primary dot) and one target circle has arrow label $2$. In this case
$\partial_{vect}\neq \partial$ because the non-trivial secondary dot
leads to either $+1$ or $-1$ in the vector grading, hence a non-trivial
higher order grading.
\bigbreak

Note that this situation does not depend on the particular choice of
chain ($1$ or $X$ in a given state). It depends only on
the labelling in the two neighboring states. We call this shifting in the vector grading the  {\em odd dotting condition}.
\bigbreak

The following Lemma follows from the definition of the vector
grading.

\begin{lm}
If for an atom we have $P=Q=R=S$ then the odd dotting condition does not
hold for any of the four edges of the bifurcation diagram.\label{lmlm2}
\end{lm}
\bigbreak

Now, it turns out that $\partial_{vect}$ in some cases can play the
role of the differential, that is, in some cases,
$\partial_{vect}^{2}=0$.
\bigbreak

Namely, we have the following
\begin{lm}
If the odd dotting condition fails, then $\partial$ does not decrease
the vector grading, that is, $\partial=\partial_{vect}+{\tilde
\partial}$ where $\partial_{vect}$ preserves the vector grading
and ${\tilde\partial}$ increases exactly one of the dotted gradings (one of the vector slots)
by $2$.\label{lmlm3}
\end{lm}
\bigbreak

\begin{proof}
We deal with a $2\to 1$-bifurcation or $1\to 2$-bifurcation. We may
assume that precisely two of the $3$ circles are dotted; moreover,
without loss of generality, we may think that these two circles have
a primary dot.
\bigbreak

Then we have to list all possible maps $m$ and $\Delta$ to see that
some of them preserve the vector grading, and the others
increase the vector grading by $2.$ Note that all calculations occur in one vector slot since the odd
dotting condition fails. In this context we can speak freely about the dotted grading and whether it increases or
decreases under a differential.
\bigbreak

Let us start with the multiplication. We see that the multiplication
of $1$ (without dot) with any of $1$, ${\dot 1}$, $X$, ${\dot X}$
leads to $1$, ${\dot 1}$, $X$, ${\dot X}$ and this multiplication
preserves the dotted grading. Now, ${\dot 1}\otimes {\dot X}$ (or
${\dot X}\otimes {\dot 1}$) multiply to get $X$, which does not
change the dotted grading. Multiplication of $X$ (or ${\dot X}$)
with another $X$ (or ${\dot X})$ gives zero. \bigbreak

Finally, ${\dot 1}\otimes {\dot 1}\to 1$ increases the dotted grading
by $2$ as well as any of ${\dot 1}\otimes X\to {\dot X}$ or
$\dot{X} \otimes 1 \to {\dot X}$.
\bigbreak

With comultiplication the situation is quite analogous. When none of
the three circles is dotted, then the dotted grading is preserved
under multiplication. If the circle in the source space and one
circle in the target space is dotted, then the comultiplication
looks like ${\dot 1}\to {\dot 1}\otimes X+{\dot X}\otimes 1$ or
${\dot X}\to {\dot X}\otimes X$. Here the only term where the dotted
grading is not preserved, is ${\dot 1}\to {\dot X}\otimes 1$; in
this case it is increased by $2$.
\bigbreak

If the circle in the source space is not dotted and both circles in
the target space are dotted then the dotting is preserved for $1\to
1\otimes X+X\otimes 1$, and it is increased by $2$ for $X\to
X\otimes X$.
\end{proof}
\bigbreak

Altogether Lemmas \ref{lmlm2} and \ref{lmlm3} lead to the following
\begin{lm}
Assume for an atom representing a face of the bifurcation cube the
labellings of all four states coincide. Then the restriction of
$(\partial')^{2}$ to this atom gives zero.
 \label{lmlm4}
\end{lm}
\bigbreak

\begin{proof}
We see that the differentials $\partial$ and $\partial_{mult}$ agree
along the edges of such an atom because of Lemma \ref{lmlm1}, so the
$2$-face corresponding to that atom $\partial_{mult}$
(anti)commutes. Moreover, by Lemma \ref{lmlm2}, the differential
$\partial$ splits into the sum of two differentials,
$\partial'+\partial''$, where $\partial''$ strictly increases the
multi-dotted grading. This means that $(\partial')^{2}=0$ because
$(\partial')^{2}$ is a composition of $(\partial^{2})=0$ with the
projection to the ``dotted-grading preserving subspace''
\end{proof}
\bigbreak

The next lemma is obvious.
\bigbreak

\begin{lm}
Assume in the setting above $P\neq S$. Then both compositions for
our atom are zero maps because of the multi-grading. Thus, the
restriction of $(\partial')^{2}$ to this atom is zero.\label{lmlm5}
\end{lm}
\bigbreak

\begin{proof}
This happens just because $\partial'$ preserves the multi-grading,
and so does $(\partial')^{2}$.
\end{proof}
\bigbreak

In the third case we have $P=Q=S\neq R$ or $P=R=S\neq Q$.
\bigbreak

In this case, we must separately consider all the six atoms
(the schema representing each atom depicted as in
Figure \ref{atm}) to show that the corresponding faces of the cube
anti-commute. We shall draw each atom separately in referring to the appropriate Figures in the paper.
\bigbreak

Consider the upper left atom depicted in Figure \ref{allatoms}.
We leave it to the reader to label the maps so that $f_{2}$
and $f_{4}$ correspond to $1\to 1$-bifurcations.
The composition $f_{4}\circ f_{2}$ is then a zero-map, by definition.
The remaining two maps are labelled $f_{1}$
and $f_{3}.$
\bigbreak

Thus, we have two options. If $R\neq P$ then the other composition
of differentials is zero because of multiple gradings. If $P=R=S$
then in the $A$-state we have only one circle labelled by $P$ as
well as in the $B$-state; in the intermediate state we have two
circles labelled by $P$ and $0$.
\bigbreak

The composition $f_{3}\circ f_{1}$ behaves as follows. First, we
comultiply $1$, and then we multiply the result. If we start with
$X$, we would end up with $0$ because $X\to X\otimes X\to 0$
even for the usual differential $\partial$. If we had $1$ then two
options are possible. If $P=0$ then the composition $f_{3}\circ
f_{1}$ will lead to $1\to 1\otimes X+X\otimes 1\to X+X=0$. If $P\neq
0$ then the $f_{3}\circ f_{1}$ will take $1$ to $0$ as well because
of the vector grading: the vector grading of $1$ for a
non-zero $P$ differs from that for $X$ by sign.
\bigbreak

Thus, for the unique non-orientable essential atom with two vertices
we have the equality $f_{4}\circ f_{2}=0=f_{3} \circ f_{1}$, which shows
the (anti)commutativity. For the other atom with the same frame (which
corresponds to the Hopf link with the $A$-state having $2$ circles)
the ``bad'' situation does not occur, just because two single-circle
states can not have different $K_{j}$'s. This completes the analysis of
the upper left atom in Figure \ref{allatoms}.
\bigbreak

We now consider the remaining five essential atoms in Figure \ref{allatoms}.
The atoms are all orientable, so the arrow count (labelling) is additive.
Following the methodology of our previous argument, we can verify that
the anticommutativity survives after the new grading is imposed for these atoms.

\begin{figure}
\centering\includegraphics[width=200pt]{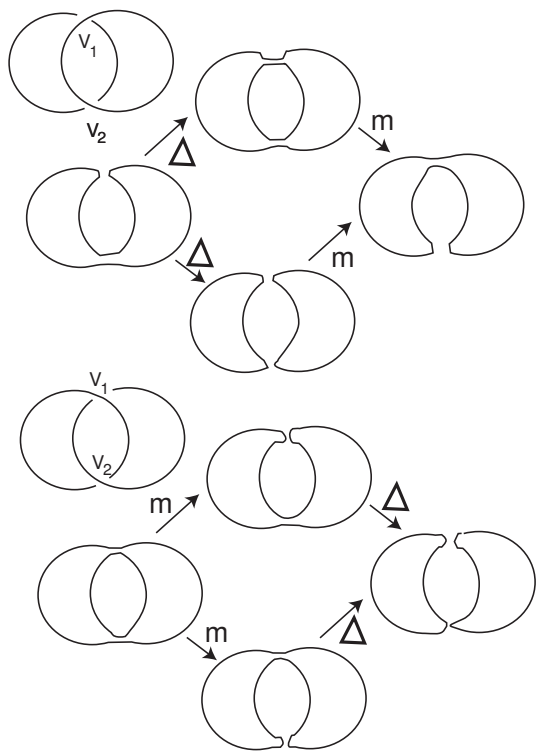}
\caption{The two atoms with the shadow of the Hopf link}

\label{Hopflink}
\end{figure}

The unlink (bottom right in Figure \ref{allatoms}) has one circle in the opposite states
and two circles in the intermediate states (see the upper part of Figure \ref{Hopflink}).
The Hopf link has 2 circles in the A-state,
2 circles in the B-state, and 1 circle in each of the two intermediate states,
as shown in the lower part of Figure \ref{Hopflink}.

\begin{figure}
\centering\includegraphics[width=250pt]{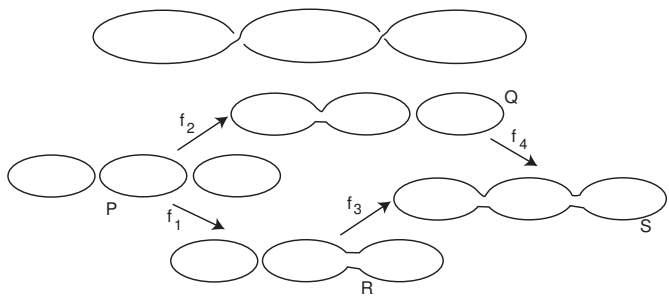}
 \caption{An atom with 2 vertices} \label{others1}
\end{figure}

Consider the three atoms having the frame of the unknot with two curls as shown in Figure \ref{allatoms}.
The corresponding bifurcation cubes have a state with three circles,
two states with two circles and one state with one circle (that is
positioned opposite the state with three circles).
The three possible bifurcation cubes depend on the number of circles in the initial
state of the cube. An example of this is shown in Figure \ref{others1}.

For the Hopf link,  assume that for both $2$-circle states
the multiple grading is the same
as that of one of the two $1$-circle states. By definition, this
means that one of the two circles in one $2$-circle state has arrow
count zero. Denote the arrow count for the other circle by $p$.
Consequently, the other way of merging the two circles gives us $p$
again. This means that the labelling is $A=B=C=D=\{p\}$, and we are
in the situation of Lemma \ref{lmlm4}.
\bigbreak

If we have $1$-circle in the $A$-state and $1$-circle in the
$B$-state, we may have a ``bad'' situation  $(P=Q=S\neq R)$ (not
covered by Lemmas \ref{lmlm4} and \ref{lmlm5})  occurring
as described below.
\bigbreak

First, note that if $P=Q=S$ then the $A$-state with two circles
should have labelling $\{P\}$ as well as the $B$-state, whence the
labelling for two circles corresponding to $Q$ should be $\{P,0\}$.
We are interested in the case when the other intermediate state has
labelling $R=(\alpha,\beta )$, say, $( \alpha,P\pm \alpha )$, where
$\alpha\neq 0$.
\bigbreak

In this case the composition $f_{3}\circ f_{1}$ is zero. Let us
consider the composition $f_{4}\circ f_{2}$.
\bigbreak

First, let us consider the partial differentials corresponding to
$\partial.$
If we apply it to $X$, we get $0$, because the
comultiplication $f_{1}$ gives us $X\otimes X$ and the further
multiplication gives zero. On the other hand, the composition
$f_{4}\circ f_{2}$ takes $1$ to $0$ because we first get $1\otimes
X+X\otimes 1$, which is then mapped to $2X=0$. Now, when we pass
from $\partial$ to $\partial'$, we see that both multiplication and
comultiplication either preserve the vector grading or
increase it by $2$, we should compare the dotting of the initial $1$
and the final $X$. If they both are zero, then the composition takes
$1$ to $2X=0$, otherwise, $1$ is taken to $0$ because of the dotted
gradings.
\bigbreak

{\it Note that this is precisely the case where we need our
coefficients to be defined over ${\bf Z}_{2}$.}
\bigbreak

Finally, all 3 atoms with the frame an unknot (drawn in the middle
of Figure \ref{allatoms}) are to be double-checked. \bigbreak

The three possibilities are: the $A$-state has $3$ circles, or it
has $2$ circles or it has $1$ circle, see Figure \ref{others1}.
\bigbreak

Assume $P=Q=S$ (the case $P=R=S$ is analogous because of the
symmetry). We claim that in this case $R=P$. Indeed, since we have
$3$ circles in the $A$-state, and one circle in the $B$-state, we
see that the labellings of the circles are $\{p,0,0\}$ in the
$A$-state and $\{p\}$ in the $B$-state. This yields that
$R=\{p,0\}=\{p\}$, and the (anti)commutativity follows from Lemma
\ref{lmlm4}.
\bigbreak

The atom when we have three circles in the $B$-state is
analogous.
\bigbreak

In fact, because of the symmetry, we can reduce these three cases to
two cases: when we have $1$ and $3$ at the ends, or when we have $2$
and $2$ at the ends.
\bigbreak

\begin{figure}
\centering\includegraphics[width=250pt]{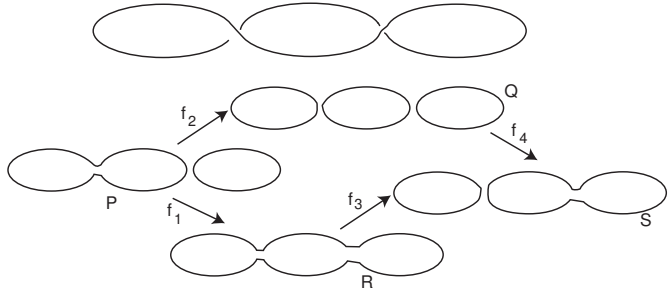}
 \caption{An atom with 2 vertices} \label{example2}
\end{figure}

Now, we are left with the example shown in Figure \ref{example2}.
\bigbreak

We are interested in the case when $P=S$ and either $Q\neq P$ or
$R\neq P$.
\bigbreak

Note that each of $P$ and $S$ consists of $2$ circles. Assume
$P=S=\{a,b\}$.
\bigbreak

It is easy to see that if $R=P=S$ then $Q=P=S=R$. Indeed, if $R=S$,
this means that both $P$ and $S$ are of the form $\{a,0\}$ (or both
are $\{b,0\}$) which yields $Q=\{a,0,0\}=P=S$.
\bigbreak

Thus we are interested in the case when $a\neq 0, b\neq 0$ and
$Q=P=S$. This means that $Q=\{a,0,b\}$, whence $R$ may be of the
form $\{|a\pm b|\}.$
In this case the composition $f_{3}\circ f_{1}=0$ because $Q\neq P$.
Let us show that the composition $f_{4}\circ f_{2}=0$.
\bigbreak

Recall that both $f_{4}$ and $f_{2}$ are compositions of the partial
differential $\partial$ with the projection map preserving the
multi-grading and the multi-dotted grading.
\bigbreak

Regardless any grading, $f_{4}\circ f_{2}$ would take $1\otimes 1\to
1\otimes X+X\otimes 1$, $1\otimes X\to X\otimes X$, $X\otimes 1\to
X\otimes X$.
\bigbreak

Now we note that none of these maps survives after applying the
projection with respect to vector grading. Indeed, consider
for instance the map from $1\otimes 1$ to the summand $1\otimes X$.
In the source space we had $1$ and $1$ with vector grading
coming from labelling $a$ and $b$; let us denote it by
$1_{a}+1_{b}$. For $1\otimes X$ we have either $1_{a}\otimes X_{b}$
or $1_{b}\otimes X_{b}$ depending on the circle having label $a$.
\bigbreak

Here $1_{a}$ denotes the $(0,\dots, 0,-1,0,\dots, 0)$ with the only
non-trivial entry on $k$-th position, $a=2^{k-1}m$ for odd $m$.
Analogously, $X_{a}$ denotes $(0,\dots,0,+1,0,\dots, 0)$ with the
only non-trivial entry on $k$-th position.
\bigbreak

It is crucially important here that neither $a$ nor $b$ is equal to
zero. This means that $1_{a}+1_{b}\neq 1_{a}+X_{b}$ just because
$1_{b}\neq X_{b}$.
\bigbreak

The same happens in the other cases.
\bigbreak

This proves that $f_{4}\circ f_{2}=0$, and the atom is
(anti)commutative because both compositions are zeroes.
\bigbreak

This completes the check of cases of the different atoms corresponding to faces of the bifurcation
cube.
\bigbreak

{\bf Part 2. Proof that the homology is invariant under Reidemeister
moves.}

Below, we shall sketch the outline of the main ideas of the proof.
The main features mirror the invariance proof for the usual
Khovanov homology along the lines of \cite{BN}.

The invariance under the first Reidemeister move is based on the
following two statements which will held when adding a small curl:

\begin{enumerate}

\item The mapping $\Delta$ is injective.

\item The mapping $m$ is surjective.

\end{enumerate}

In fact, the last two conditions hold when the small circle has the trivial
arrow count, and this means that it does not contribute to any of the gradings.
\bigbreak

Indeed, consider the complex

\begin{equation} [[\skkinkr]]=\left([[\skroh]]\stackrel{\Delta}{\to} [[\skrov]] \{1\}\right).\end{equation}

The usual argument goes as follows: the complex in the right hand
side contains a $\Delta$-type partial differential, which is
injective. Thus, the complex $[[\skrov]]$ is killed, and what remains
from $[[\skroh]]$ is precisely (after a suitable normalisation) the
homology of $[[\skkinkl]]$.
\bigbreak

But $\Delta$ is injective because for any $l\in \lbrace 1,X \rbrace$ we have
$\Delta(l)=l\otimes X+ \langle \mbox{other terms} \rangle$, where the second
term $X$ in $l\otimes X$ corresponds to the small circle.
\bigbreak

But in our situation with dotted circles, this happens only if {\em
the small circle is not dotted}. But if the small circle has
non-trivial arrow count (say, it appears after splitting a circle
without dots into two circles with primary dot each), it would lead,
say, to $\Delta:X\to 0$, because $\mbox{$\dot X$}\wedge \mbox{$\dot X$}$ has another
vector grading (which is greater by $2$ than the grading of $X$).
\bigbreak

An analogous situation happens with the other curl
\begin{equation} [[\skkinkl]]=\left([[\skrov]]\stackrel{m}{\to} [[\skroh]] \{1\}\right).\end{equation}

Here we need that  the mapping $m$ be surjective; actually, it would
suffice that the multiplication by $1$ on the small circle is the
identity. But this happens if and only if the small circle has arrow count
0, that is, we have $1$, not $\dot{1}$.
\bigbreak

Quite similar things happen for the second and for the third
Reidemeister moves. The necessary conditions can be summarised as
follows:
\bigbreak

{\em The small circles which appear for the second and the third
Reidemeister move should not be dotted}, and similar  curves for corresponding  smoothings of the RHS and the LHS
of any Reidemeister move have the same dotting.
\bigbreak

The explanation comes a bit later.  Now, we see that this condition
is obviously satisfied when the dotting comes from a cohomology
class, and not necessarily the Stiefel-Whitney cohomology class for
non-orientable surface. Any homology class should do.
\bigbreak

Thus (modulo some explanations given below) we have proved the
following
\bigbreak

\begin{thm}
Let ${\cal M}\to M$ be a fibration with $I$-fibre so that ${\cal M}$
is orientable and $M$ is a $2$-surface. Let $h$ be a
$Z_{2}$-cohomology class and let $g$ be the corresponding dotting.
Consider the corresponding grading on $[[K]]$. Then for a link
$K\subset {\cal M}$ the homology of $[[K]]_{g}$ is invariant under
isotopy of $K$ in $M$ (with both the orientation of $M$ and the
$I$-bundle structure fixed) up to some shifts of the usual (quantum)
grading and height (homological grading).
\end{thm}

{\bf Explanation for the second and the third moves.}

We have the following picture for the Reidemeister move for
$[[\skrtwhh]]$:
\bigbreak

\begin{equation}\begin{array}{ccc}
\lb\skrtwhh\rb\{1\} & \stackrel{m}{\longrightarrow} & \lb\skrtwvh\rb\{2\}\\
\Delta\uparrow &  & \uparrow \\
\lb\skrtwhv\rb & \longrightarrow & \lb\skrtwvv\rb\{1\}
\end{array}.
\end{equation}

Here we use the notation ${\{\cdot \}}$ for the degree shifts, see
page \ref{shifts}.
\bigbreak

\begin{equation}\begin{array}{ccc}
\lb\skrtwhh\rb\{1\} & \stackrel{m}{\longrightarrow} & \lb\skrtwvh\rb\{2\}\\
\Delta\uparrow &  &\uparrow \\
\lb\skrtwhv\rb & \longrightarrow & \lb\skrtwvv\rb\{1\}
\end{array}.
\end{equation}
This complex contains the subcomplex ${\cal C}'$:

\begin{equation} {\cal C}'=\begin{array}{ccc}
\lb\skrtwhh\rb_{1}\{1\} & \stackrel{m}{\longrightarrow} &\lb\skrtwvh\rb\{2\}\\
 \uparrow & & \uparrow \\
0 & \longrightarrow & 0
\end{array}\end{equation}
if {\em the small circle is not dotted}.
\bigbreak

From now on $1$ denotes the mark on the small circle.
Then the acyclicity of ${\cal C}'$ is evident.
Factoring ${\cal C}$ by ${\cal C'}$, we get:
\bigbreak

\begin{equation}\begin{array}{ccc} \lb\skrtwhh\rb\{1\}\slash_{1=0} &
\longrightarrow&  0\\
\Delta\uparrow & & \uparrow \\
\lb\skrtwhv\rb & \longrightarrow &
\lb\skrtwvv\rb\{1\}\end{array}.\end{equation}

In the last complex, the mapping $\Delta$ directed upwards, is an
isomorphism (when our small circle is not dotted). Thus the initial
complex has the same homology group as $[[\skrtwvv]]$. This proves
the invariance under $\Omega_{2}$.
\bigbreak

The argument for $\Omega_{3}$ is standard as well; it relies on the
invariance under $\Omega_{2}$ and thus we also require that
the small circle is not dotted.
\bigbreak

\section{Applications}

The complex constructed in this paper allows us to prove some properties of virtual
knot diagrams coming from the Kauffman bracket, the Khovanov
homology and the arrow polynomial, see \cite{MyBook},\cite{ExtBr},
\cite{Minimal},\cite{TuraevGenus},\cite{DK}.
\bigbreak

First, the consideration of the chain spaces and arrow counts
immediately leads to the following theorem.
\begin{thm}
Assume K is a virtual link diagram, and assume there is a
non-trivial homology class of [[K]] with multiple grading
$\{k_{1},\dots, k_{n}\}$, such that $\sum_{i} |k_{j}|=k$. Then any
diagram of $K$ has at least $k$ virtual crossings.
\end{thm}

Besides, the following generalization of the Kauffman-Murasugi
Theorem says
\begin{thm}
Let $K$ be a virtual link diagram with a connected shadow (that is,
every classical crossing of $K$ can be connected to any other
classical crossing by a sequence of arcs starting and ending at
classical crossings and going through virtual crossings).
\smallbreak

Let $g$ be the minimal  oriented atom genus for the diagram of $K$ and
let $n$ be the number of crossings in the diagram $K$.
Then $span \langle K \rangle\le 4n-4g$, where $span$ stands for the
difference between the leading degree and the lowest degree of the
Kauffman bracket with respect to the variable $a$.\label{gnrl}
\end{thm}

The condition of theorem \ref{gnrl} rules out the split link
diagrams.
The same argument (see \cite{MyBook,Minimal,Doklady}) leads to
\begin{thm}
For $K$ as in Theorem \ref{gnrl}, the span of the arrow polynomial
of $K$ taken with respect to $a$  does not exceed
$4n-4g$.\label{thm}
\end{thm}

On the other hand, the genus of the atom estimates from above the
{\em thickness of the Khovanov homology}: the number of diagonals
with slope two in coordinates (homological grading, quantum grading)
which appear between the leftmost and the rightmost diagonal having
a non-trivial homology group. The estimate in \cite{Doklady} says
that this thickness does not exceed $2+g$.
Similar considerations lead to the same estimate for the thickness
of $[[\cdot]]$ (taking with respect to the {\em old gradings}, after
forgetting all new gradings of non-trivial homology groups):
\bigbreak

\begin{thm}
For $K$ as in Theorem \ref{gnrl}, the thickness of $[[\cdot]]$ does
not exceed $2+g$.\label{thj}
\end{thm}

Theorems \ref{thm} and \ref{thj} together lead to the following
\begin{thm}
Assume the diagram $K$ represents a split virtual link (e.g. virtual
knot). Then,  if $K$ having $span$ of the arrow bracket equal to
$4n-4g$ and the thickness of the extended Khovanov homology equal to
$2+g$ then this diagram is minimal with respect to the number of
classical crossings.
\end{thm}

It is an interesting question to
determine if there exist examples where the theorems stated above
give sharper estimates than the already existing invariants.
\bigbreak

\section{Open questions}

The methods described in the present paper allow us to extend the
arrow counts in the arrow polynomial to the level of gradings
of a link homology theory. We can recover the arrow polynomial from this link
homology by taking the Euler characteristic, forgetting vector gradings
and taking the multiple gradings as arrow counts. In this sense, our link
homology theory is a true categorification of the arrow polynomial.
\bigbreak

There is a more delicate invariant, the {\em extended bracket
polynomial}, \cite{ExtBr}, which generalizes the arrow polynomial
and takes geometrical information into account (instead of just
arrow counts).
Can this polynomial be categorified by using techniques given in the
present paper?
\bigbreak

Another question is whether there is a categorification of the arrow
polynomial (or the extended bracket polynomial) with integral
coefficients. The only point where we needed the ${\bf Z}_{2}$
coefficients was the atom in Fig. \ref{others1} where the vector
gradings and the multiple gradings together did not make the complex
over ${\bf Z}$ well defined.
However, in a similar situation one gets the commutativity of the
corresponding face of the atom for {\em odd Khovanov homology
theory}, \cite{ORS}.
Thus, the question of generalizing odd Khovanov homology theory for
virtual links gets one more motivation: it would be useful to have it for constructing a
categorification of the arrow polynomial with integral coefficients.
\bigbreak

Another issue of investigation is the notion of {\em parity of
crossings}, developed recently by Manturov, \cite{Parity,Parity2}
(see also \cite{KaC} for a precursor to this approach).
The idea is to distinguish between two types of
crossings, the {\em even ones}, and the {\em odd ones} according to
some axioms. This approach turns out to be extremely powerful in
recognizing some virtual knots and creating new virtual knot
invariants. There is a natural way to generalize the arrow
polynomial by using the parity argument. This, and a corresponding categorification will be discussed
in a subsequent paper.
\bigbreak

\end{document}